\documentclass[a4paper,11pt]{amsart}

\usepackage{txfonts}
\usepackage[T1]{fontenc}
\usepackage{fancyhdr}

\usepackage{amsthm}
\usepackage{amsmath}
\usepackage{enumerate}
\usepackage{amsfonts}
\usepackage{amssymb}
\usepackage{fullpage}
\usepackage{amsmath,amscd}
\usepackage{stmaryrd}
\usepackage{graphicx}
\usepackage{tikz-cd}
\usepackage{array}
\usepackage{amsmath,amssymb,amsfonts,dsfont}
\usepackage[english]{babel}
\usepackage[T1]{fontenc}
\usepackage{graphicx}
\usepackage{float}
\usepackage{pdfpages}
\usepackage[a4paper, margin = 4cm, bottom = 3.5cm]{geometry}

\usepackage{ifpdf}
\usepackage{marginnote}
\usepackage{hyperref}	

\usetikzlibrary{calc}
\usetikzlibrary{decorations.pathreplacing,decorations.markings,decorations.pathmorphing}
\usetikzlibrary{positioning,arrows,patterns}
\usetikzlibrary{cd}
\usetikzlibrary{intersections}
\usetikzlibrary{arrows}

\usepackage{color}
\usepackage{hyperref}
\hypersetup{
    colorlinks=true, 
    linktoc=all,     
    linkcolor=blue,  
    citecolor=blue,
    filecolor=blue,
    urlcolor=blue
}

\newcommand{\supr}[1]{\underset{#1}{\sup}}

\newcommand{\CC}{{\mathbb{C}}}
\newcommand{\RR}{{\mathbb{R}}}
\newcommand{\NN}{\mathbb{N}}

\newtheorem{thm}{Theorem}[section]

\newtheorem{cor}[thm]{Corollary}
\newtheorem{prop}[thm]{Proposition}
\newtheorem{lem}[thm]{Lemma}

\theoremstyle{definition}
\newtheorem{defn}[thm]{Definition}
\theoremstyle{remark}
\newtheorem{rem}[thm]{Remark}
\theoremstyle{definition}

\theoremstyle{definition}

\theoremstyle{definition}
\newtheorem{problem}[thm]{Problem}
\theoremstyle{definition}

\theoremstyle{definition}

\numberwithin{equation}{section}

\title[Closed geodesics in dilation surfaces]{Closed geodesics in dilation surfaces}

\author{Adrien Boulanger}
\address[Adrien Boulanger]{Institut mathématique de Marseille, Aix-Marseille Université, Technopôle Château-Gombert 39, rue Frédéric Joliot-Curie, Marseille}
\email{adrien.boulanger@univ-amu.fr}

\author{Selim Ghazouani}
\address[Selim Ghazouani]{Department of Mathematics, University College London, United Kingdom}
\email{s.ghazouani@ucl.ac.uk}

\author{Guillaume Tahar}
\address[Guillaume Tahar]{Beijing Institute of Mathematical Sciences and Applications, Huairou District, Beijing, China}
\email{guillaume.tahar@bimsa.cn}

\date{\today}

\begin{document}

\begin{abstract}
We prove that directions of closed geodesics in every dilation surface form a dense subset of the circle. The proof draws on a study of the degenerations of the Delaunay triangulation of dilation surfaces under the action of Teichm\"{u}ller flow in the moduli space.
\end{abstract}
\maketitle
\setcounter{tocdepth}{1}
\tableofcontents

\section{Introduction}

We consider in this article the problem of the existence of regular closed geodesics in dilation surfaces. Our main theorem is the following.

\begin{thm}
\label{thm:main}
For any closed dilation surface $\Sigma$, there is a dense set of directions $\theta$ such that the directional foliation $\mathcal{F}_{\theta}$ has a periodic leaf. Equivalently, the set of directions covered by a cylinder is dense in $\mathbb{RP}^{1}$.
\end{thm}

In particular, any dilation surface $\Sigma$ carries at least one closed geodesic. This generalizes to the context of dilation surfaces a celebrated theorem of Masur \cite{Mas86} for translation surfaces.

\vspace{3mm} As the two equivalent formulations of Theorem \ref{thm:main} suggest, it can be viewed from either a dynamical or geometric perspective. From the geometric point of view, it guarantees that every dilation surface contains the simplest building block that can be imagined, a cylinder, thus giving valuable insight into the geometric structure of the arbitrary dilation surface.
 
\vspace{2mm} \noindent On the dynamical side, this theorem guarantees the ubiquity of periodic orbits in some particular (but very natural) one-parameter families of unidimensional of dynamical systems, in the form of the following corollary.

\begin{cor}
\label{cor:AIT}
For every affine interval exchange transformation $T_0 : [0,1] \longrightarrow [0,1]$, the set of parameters $t$ such that the map $x \mapsto T_0(x) + t \ \mathrm{mod} \ 1$ has a periodic orbit is dense in $\mathbb{R}$.
\end{cor}

Results about particular families of dynamical systems of this type are usually difficult to prove; a result analogous to Corollary \ref{cor:AIT} where $T_0$ is an arbitrary generalized interval exchange map seems out of reach of current methods. 

\subsection{Affine structures on surfaces} The question of the existence of closed geodesics can be considered in the wider context of affine (complex or real) structures on surfaces\footnote{~A real (resp. complex) affine structure on a surface is an atlas of charts taking values in $\mathbb{R}^2$ (resp. $\mathbb{C}$) such that transition maps lie in the group of real affine transformations $GL^{+}(2,\mathbb{R}) \ltimes \mathbb{R}^2$ (resp. complex affine transformations $\mathbb{C}^{\ast} \ltimes \mathbb{C}$), with possibly finitely many cone-type singularities.}. For Riemannian structures, the existence of closed geodesics has been known for a long time (see for example \cite{GM69} and references therein). The case of translation surfaces, which lies in the intersection of the affine and Riemannian world, is now very well understood. On the contrary, for general affine structures very little is known. We therefore pose the following problem.

\begin{problem}
Characterise the affine structures on closed surfaces which carry a regular\footnote{~A regular closed geodesic is a closed geodesic that does not contain any singularity of the affine structure.} closed geodesic.
\end{problem}

Note that a complete solution to this problem is likely to be very difficult, as it contains as a particular case the notoriously hard question of determining whether the billiard flow of every polygonal table has a periodic orbit.



\subsection{Dilation surfaces vs general affine surfaces} Dilation surfaces are particular complex affine surfaces whose structural group is the set of transformations of the form $ z \mapsto a \cdot z + b$ where $a$ is a positive real number and $b \in \mathbb{C}$. Although it is expected that generic complex affine surfaces do not have any closed geodesic, the main theorem we prove in this article predicts that any dilation surface does.\newline

We explain what the condition on the structural group defining dilation surfaces implies at the dynamical level. Every (complex or real) affine structure induces a geodesic foliation on $\mathrm{T}^1 \Sigma$ the unit tangent bundle of the surface. $\mathrm{T}^1 \Sigma$ is a tridimensional manifold thus the dynamical system induced by the foliation is essentially bidimensional. Indeed, for a given Poincar\'{e} section, the first return map may change both the direction and the position of the intersection of the leaf with the interval.\newline

In the particular case of dilation surfaces, $\mathrm{T}^1 \Sigma$ decomposes into a $1$-parameter family of invariant surfaces for the foliation. While this gives no indication as to which affine structures always have periodic leaves, it explains why dilation surfaces are essentially different from the general case: 
\begin{itemize}
\item the problem for dilation surfaces is about finding periodic orbits in \textit{one-parameter families} of unidimensional dynamical systems;
\item the problem for the generic affine surface is about finding a periodic orbit for a \textit{given} bidimensional dynamical system.
\end{itemize} 

The analysis of bidimensional dynamical systems is far more intricate than that of their unidimensional counterparts; furthermore with dilation surfaces we have an entire one-parameter family of unidimensional dynamical systems (which are easier to analyse) to find a periodic orbit. This discussion also explains why, despite the fact that in principle it is plausible that a lot of real affine surfaces carry closed geodesics, the dilation case is of a different nature and probably easier to analyse.

\subsection{Action of $SL(2,\mathbb{R})$ and strategy of proof.}

We now explain the ideas behind the proof of Theorem \ref{thm:main}. It is very much inspired by the translation case, and we remind the reader of the general structure of the proof in this case. We refer to \cite{Mas86} for the original proof in the translation case. 

\vspace{2mm} Both moduli spaces of dilation and translation surfaces carry an action of the group $SL(2,\mathbb{R})$. This action is naturally defined by the post-composition of the charts defining the dilation/translation structure. It has the following important property: two surfaces are on the same $SL(2,\mathbb{R})$-orbit if and only if they define the same underlying real affine structure. In particular, if a surface has a closed geodesic, it is the case for every surface in its $SL(2,\mathbb{R})$-orbit.

\vspace{2mm} In the translation case, the proof goes by induction on the combinatorial complexity of the surface.\footnote{~We define the complexity to be the number of triangles in a triangulation whose set of vertices is the set of singular points of the surface.}

\begin{enumerate}
 \item It is easy to check that translation surfaces of lowest complexity (flat tori) always carry closed geodesics. 

\item Assume that we know that all surfaces of complexity lesser than $k$ do carry closed geodesics and consider a translation surface $\Sigma$ of complexity $k$. It is not hard to find a sequence $(\Sigma_{n})_{n \in \mathbb{N}}$ of translation surfaces in the $SL(2,\mathbb{R})$-orbit of $\Sigma$ which diverges \textit{i.e.} leaves any compact subset in the moduli space of surfaces of complexity $k$. 

\item Geometric tools building on the Riemannian structure of translation surfaces allow to show the following dichotomy: either $(\Sigma_n)_{n\in\mathbb{N}}$ Gromov-Hausdorff converges (up to passing to a subsequence) towards a translation surface of lesser complexity or the Riemannian diameter of $\Sigma_{n}$ tends to infinity.

\item In the first case, having a cylinder is a property that is open in parameter space and by the induction hypothesis, for $n$ large enough $\Sigma_{n}$ has a closed geodesic. Since $\Sigma_{n}$ has the same real affine structure as $\Sigma$, so does $\Sigma$. 

\item In the second case, an elegant lemma due to Masur and Smillie (Corollary 5.5 in \cite{MS91}) ensures that a translation surface of large diameter contains a long flat cylinder and thus contains a closed geodesic, which concludes the proof.

\end{enumerate}

This strategy relies heavily on the Riemannian nature of translation surfaces to get a rather simple analysis of the ways a sequence of translation surfaces can degenerate; this part of the proof breaks down when trying to generalize it to the case of dilation surfaces. Most of the work done in the present article is to replace the last three points of the strategy outlined above by a suitable analysis of the different ways a sequence of dilation surfaces can degenerate. We will give a precise roadmap of the proof in Section~\ref{sec:overview}. The three key technical steps of the proof (Propositions~\ref{prop.degenerating},~\ref{prop.maximaldomains} and~\ref{prop.petiteportegrandeporte}) are proved respectively in Section~\ref{sec:divergent},~\ref{sec:control} and~\ref{sec.petiteportegrandeporte}.

\subsection{An important shortcoming and an open problem}

We prove that every dilation surface contains a closed geodesic but unfortunately we were not able to infer anything concerning the nature of the cylinder carrying this closed geodesic. In particular, our proof does not preclude the existence of a dilation surface which is not a translation surface all of whose cylinders are flat (although the existence of such a surface seems highly unlikely). 

\begin{problem}
Show that a dilation surface whose cylinders are all flat is a translation surface.
\end{problem}

\paragraph{\bf Acknowledgements.} The second author is greatly indebted to Bertrand Deroin for introducing him to the topic of affine structures on surfaces and asking him the question that lead to the present article. The third author would like to thank Dmitry Novikov for interesting feedback. The authors are grateful to the anonymous referee for valuable remarks and discussions.

\section{Dilation surfaces}\label{sec:dilation}

The symbol $\Sigma$ will always stand for a compact surface of genus $g \geq 0$ with a finite number of boundary components. 

\subsection{Dilation cones} Singularities of dilation surfaces are modelled on singularities of dilation cones.

\vspace{3mm} For any $k \in \mathbb{N}^{\ast}$, a \textbf{flat cone} of angle $2k\pi$ is obtained as the cyclic cover of $\mathbb{C}$ of degree $k$ ramified at $0$. Vertex $0$ is a \textbf{cone point} of angle $2k\pi$ in the flat cone.

\vspace{3mm} For any $k \in \mathbb{N}^{\ast}$ and any $\lambda \in \mathbb{R}^{\ast}$, a \textbf{dilation cone} of angle $2k\pi$ and multiplier $\lambda$ is obtained from a flat cone of angle $2k\pi$ by cutting a slit along a half-line starting from vertex $0$ and identifying the left side with the right side by a homothety of multiplier $\lambda$. Vertex $0$ is then a \textbf{cone point} of angle $2k\pi$ and dilation multiplier $\lambda$.

\vspace{2mm}In particular, for the affine structure induced by the gluing, the holonomy of any closed simple loop around the vertex is a homothety of dilation multiplier $\lambda$.

\subsection{Generalities} The main objects we will deal with in this article are dilation structures, defined as follows.

\begin{defn}
	\label{defn:affinesurface}
A \textbf{marked topological surface} is a topological surface $\Sigma$ - possibly with boundary - with a non-empty finite set $S \subset \Sigma$ of \textbf{marked points} such that each boundary component contains an element of $S$.

A \textbf{dilation structure} on a marked topological surface $(\Sigma,S)$ is an atlas of charts $\mathcal{A} = (U_i, \varphi_i)_{i \in I}$ on $\Sigma \setminus S$ such that
	\begin{itemize}
\item the transition maps are locally restriction of elements of
    	$\mathrm{Aff}_{\RR^*_+}(\CC) = \{ z \mapsto az+b \ | \ a \in \RR^
    *_+ \ , \ b \in \mathbb{C} \} $;
\item each marked point in the interior of $\Sigma$ has a punctured neighbourhood which is affinely equivalent to a punctured neighbourhood of the cone point
of a dilation cone;
\item each marked point on the boundary of $\Sigma$ has a punctured neighbourhood which is affinely equivalent to a neighbourhood of the centre of an Euclidean angular sector of arbitrary angle;
\item unless it is a marked point, each point of the boundary of $\Sigma$ has a punctured neighbourhood which is affinely equivalent to a neighbourhood of the centre of an Euclidean angular sector of angle $\pi$.
	\end{itemize}
Elements of $S$ are the \textbf{singularities} of the dilation structure.
\end{defn}

\begin{figure}
\includegraphics[scale=0.4]{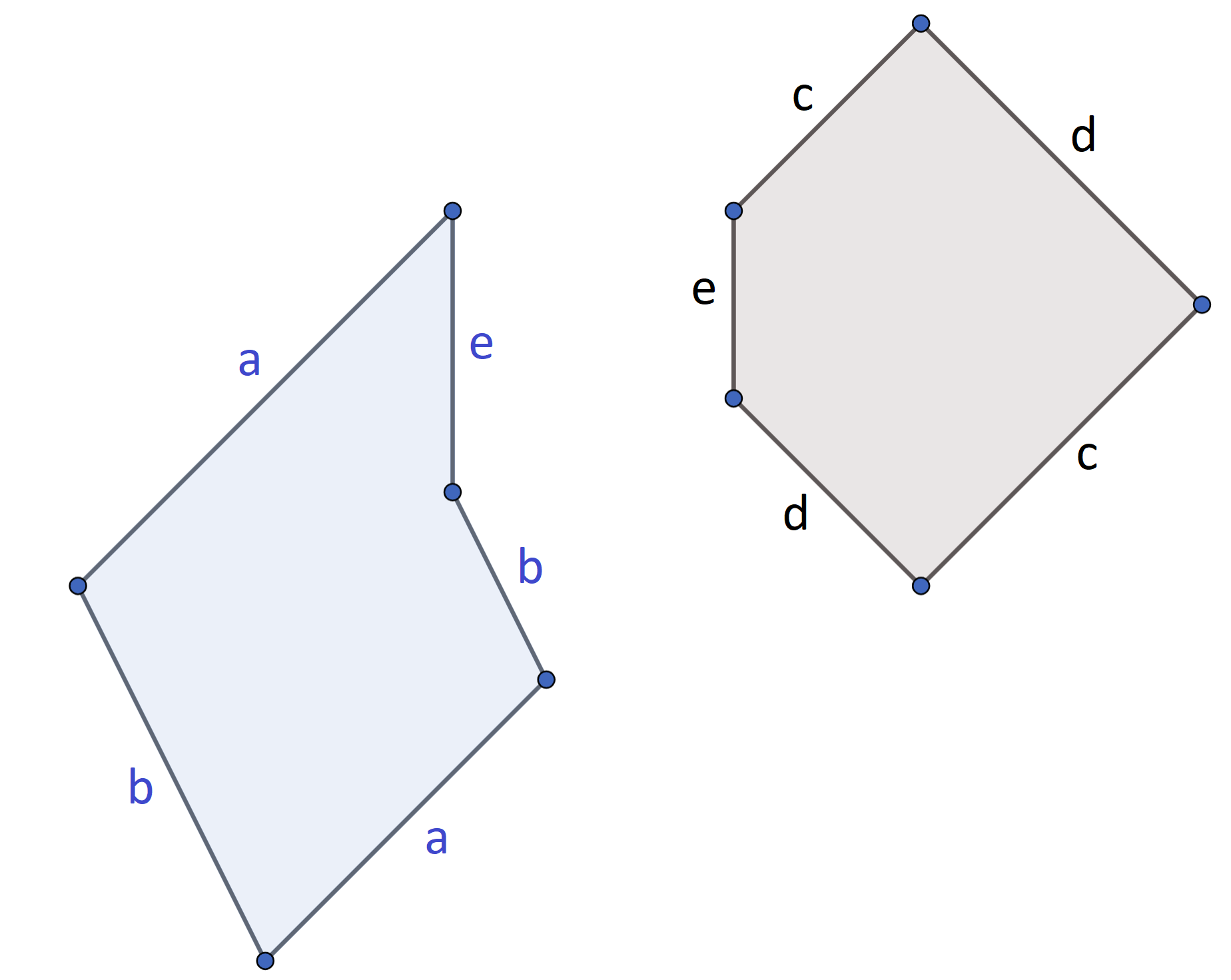}
\caption{The sides of the two polygons are glued according to their names. Topologically, the resulting surface has genus two and has only one singularity which corresponds to the extremal points of these two polygons.}
\label{figure1}
\end{figure}

A particularly simple way of constructing a dilation surface is to glue planar polygons together by using translations and dilations as illustrated in Figure~\ref{figure1}. We will see that, up to addition of finitely many singularities with an angle of $2\pi$ and a trivial dilation multiplier, every dilation surface can be constructed in this way.

Note that the notion of straight line on the surface is well-defined since changes of coordinates are affine maps. Moreover, in any direction $\theta \in \mathbb{RP}^1$ the foliation by straight lines of $\mathbb{C}$ in the direction defined by $\theta$ being invariant by dilation maps, it gives rise to a well-defined oriented foliation $\mathcal{F}_{\theta}$ on any dilation surface. Such a foliation is called a \textbf{directional foliation}. We call the \textbf{directional foliations} the resulting family of foliations, therefore indexed by $\mathbb{RP}^{1}$, that we denote $(\mathcal{F}_{\theta})_{\theta \in \mathbb{RP}^1}$. We shall call a \textbf{trajectory} any oriented leaf of one of these foliations. 
\begin{defn}

Let $\Sigma$ be a dilation surface.

\begin{itemize}
\item A \textbf{closed geodesic} in $\Sigma$ is a periodic leaf of a directional foliation. 
\item A \textbf{saddle connection} is a topological segment on the surface $\Sigma \setminus S$ which is also a straight line (a piece of leaf of a directional foliation) and whose boundary consists of two singularities (possibly identical).
\end{itemize}
\end{defn}

We conclude this subsection by the following definition we will use to measure the complexity of a dilation surface. 

\begin{lem}
We consider a compact topological surface $X$ of genus $g$ with $b$ boundary components, $n_{i}$ marked points in its interior and $n_{b}$ marked points on its boundary.\newline
Assuming $n_{b}+n_{i} \geq 1$ and that every boundary component contains at least one marked point, any topological triangulation of $X$ whose set of vertices coincides with the marked points of $X$ is formed by exactly $4g+2n_{i}+2b+n_{b}-4$ topological triangles.
\end{lem}

\begin{proof}
The Euler characteristic $\chi(X)$ of surface $X$ is $2-2g-b$. For any such topological triangulation, the number of vertices is $n_{i}+n_{b}$. Thus, we have $2-2g-b=T-A+n_{i}+n_{b}$ where $T$ is the number of triangles in the triangulation and $A$ is the number of arcs.\newline
Connected components of the boundary are loops. Thus, the number of boundary arcs is exactly $n_{b}$. Every arc has two sides (excepted the boundary arcs). Thus, we have $3T=2A-n_{b}$. We have $4-4g-2b=2T-2A+2n_{i}+2n_{b}$. It follows that $4-4g-2b=-T+2n_{i}+n_{b}$ and thus $T=4g+2n_{i}+2b+n_{b}-4$.
\end{proof}

\begin{defn}\label{defn:complexity}
The \textbf{complexity} of a marked topological surface is the number of triangles of any topological triangulation whose set of vertices is exactly the set of marked points. By convention, we define the complexity of the empty set to be zero.

We define the complexity of a dilation surface as the complexity of the underlying marked topological surface.
\end{defn}

\subsection{Cylinders}
\label{sec:cylinders} Cylinders are the geometric counterpart of the periodic leaves of the directional foliations as, in particular, each cylinder contains a closed geodesic. Conversely, any neighborhood of a closed geodesic contains a portion of cylinder. Note that we always, in this article, understand cylinders as maximal: we say that a cylinder is maximal if it is not included in any cylinder but itself.

\vspace{2mm} A \textbf{flat cylinder} is a dilation surface with boundary obtained by gluing a pair of the opposite sides of a parallelogram embedded in $\mathbb{R}^2$.

\begin{figure}
\includegraphics[scale=0.6]{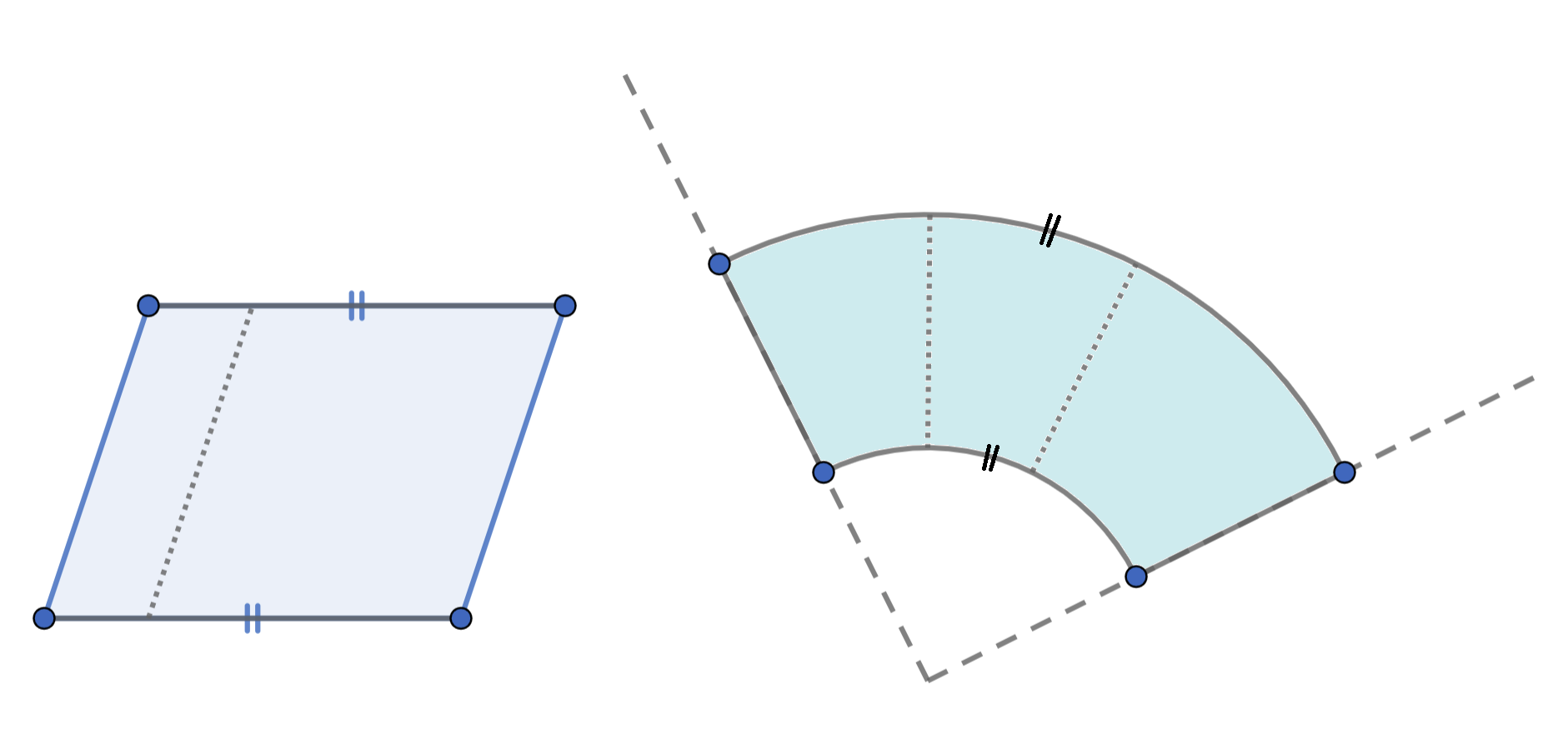}
\caption{On the left a flat cylinder with a closed geodesic in dashed (corresponding to the only direction in which there is a closed geodesic). On the right a dilation cylinder with two closed geodesics of two different directions.}
\label{figure2}
\end{figure}

\vspace{2mm} A \textbf{dilation cylinder} is a dilation surface (with boundary) obtained by cutting a sector $\mathcal{C}_{\theta}$ of angle $\theta$ in the universal cover of $\mathbb{C}^{\ast}$. The quotient of $\mathcal{C}_{\theta}$ by the dilation $z \mapsto \lambda z$ with $\lambda > 1$ real is called a \textbf{dilation cylinder} (see Figure~\ref{figure2}). 

\subsection{Moduli of cylinders}

In this paragraph we give an interpretation of conformal moduli of cylinders in dilation structures. 

\begin{itemize}

\item Recall that the modulus of a flat cylinder obtained from a rectangle of base $(z_{1},z_{2}) \in \mathbb{C}^{2}$ where the sides glued together are those corresponding to $z_2$ is by definition $\frac{|z_{2}|}{|z_{1}|}$.

\item A dilation cylinder of angle $\theta$ and dilation multiplier $\lambda > 1$ is biholomorphic (using the exponential map) to the flat cylinder obtained from a rectangle of base $(ln(\lambda),i \theta)$. Its conformal modulus is thus $\frac{\theta}{ln(\lambda)}$.

\end{itemize}

We call a closed geodesic within a cylinder a waist curve of this cylinder. 

\begin{lem}
\label{lem:intersect}
There is an absolute constant $M>0$ such that for any pair of cylinders $C_{1}, C_{2}$ of conformal modulus at least $M$ in a dilation surface $\Sigma$, either $C_{1}$ and $C_{2}$ are disjoint or their waist curves are in the same homotopy class.
\end{lem}

\begin{proof}
A consequence of Margulis lemma is the existence of an universal constant $\epsilon$ such that in any hyperbolic surface of unit area, closed geodesics of length smaller than $\epsilon$ are automatically disjoint (see Section 4.2.4 in \cite{Martelli} for a reference).\newline
In the conformal class of $\Sigma$ (punctured at the singularities), we consider the unique hyperbolic metric of unit area. In the homotopy class of waist curves of cylinder $C_{1}$ (resp. $C_{2}$), there is a unique simple closed geodesic $\gamma_{1}$ (resp. $\gamma_{2}$). We denote by $l(\gamma_{1})$ and $l(\gamma_{2})$ the lengths of geodesics $\gamma_{1}$ and $\gamma_{2}$. Assuming that waist curves of $C_{1}$ and $C_{2}$ do not belong to the same homotopy class, $\gamma_{1}$ and $\gamma_{2}$ are distinct.\newline
Following the interpretation of conformal modulus in terms of extremal length, if the conformal modulus of $C_{1}$ is strictly bigger than $M=\epsilon^{-2}$, then $l(\gamma_{1})<\epsilon$. The same holds for $C_{2}$ and $l(\gamma_{2})$. Thus, $\gamma_{1}$ and $\gamma_{2}$ are disjoint and waist curves of $C_{1}$ and $C_{2}$ do not intersect.
\end{proof}

\begin{cor}\label{cor:moduli}
For any dilation surface $\Sigma$, there exists a constant $M(\Sigma)>0$ such that any cylinder in $\Sigma$ has conformal modulus smaller than $M(\Sigma)$.
\end{cor}

\begin{proof}
We assume by contradiction that $\Sigma$ contains an infinite family of cylinders of arbitrary large moduli. Since two different cylinders always define two different free homotopy classes, we can always find intersecting cylinders with arbitrarily large moduli. This contradicts Lemma~\ref{lem:intersect}.
\end{proof}

\subsection{Pencils}\label{sec:pencil}
They have already been studied in \cite{Ta23} to make explicit the geometric properties of strict dilation surfaces in comparison with translation surfaces. We gather the needed results of this article in this section. 

\begin{defn}
A \textbf{pencil} is a continuous family of oriented trajectories starting from the same point. Let $x$ be a (possibly singular) point of a dilation surface $\Sigma$ be a dilation surface and $I$ an open interval of $\mathbb{RP}^{1}$. The notation $P(x,I)$ will refer to a pencil of trajectories starting at $x$ and covering directions of $I$.\newline
It should be noted that there are usually several pencils for a given pair $(x,I)$.
\end{defn}

The following statement, proven as Lemma~3.3 in \cite{Ta23}, provides a geometric criterion for the existence of dilation cylinders.

\begin{lem}\label{lem:cone}
Let $x$ be a point in a dilation surface $\Sigma$ (possibly with boundary) and $I$ be an open interval of $\mathbb{RP}^{1}$. For a given pencil $P(x,I)$, at least one of the following statements must hold:
\begin{enumerate}
   \item a trajectory of $P(x,I)$ hits a singularity;
    \item there exists a closed geodesic whose direction belongs to interval $I$;
    \item there is an open subset $J \subset I$ such that trajectories of the  restricted pencil $P(x,J)$ cross the interior of a boundary component of $\Sigma$.
\end{enumerate}
\end{lem}

Note that in the case where $\Sigma$ is without boundary then only the two first items can hold. \\

A second result, proven as Corollary~4.6 in \cite{Ta23}, proves the existence of dilation cylinders in dilation surfaces with non-empty boundary where trajectories of a pencil avoid the boundary.

\begin{prop}\label{prop:boundaryexceptional}
Let $\Sigma$ be a connected dilation surface with a non-empty boundary, a point $x \in \Sigma$ and an open interval $I$ in $\mathbb{RP}^{1}$. Then at least one of following statements holds 
\begin{enumerate}
    \item there is an open interval $J \subset I$ such that every trajectory of the restricted pencil $P(x,J)$ accumulates on a closed geodesic of a dilation cylinder of $\Sigma$;
    \item there is an open interval $J \subset I$ such that every trajectory of $P(x,J)$ crosses the interior of a boundary saddle connection of $\Sigma$.
\end{enumerate} 
\end{prop}

\subsection{Non-polygonable surfaces}

\begin{defn}
A polygonation of a dilation surface $\Sigma$ is family of saddle connections $\gamma_{1},\dots,\gamma_{k}$ such that:
\begin{itemize}
    \item[(i)] interiors of saddle connections $\gamma_{1},\dots,\gamma_{k}$ are disjoint;
    \item[(ii)] connected components of $\Sigma \setminus \bigcup\limits_{i=1}^{k} \gamma_{i}$ are flat polygons without any interior singularity.
\end{itemize}
A surface $\Sigma$ is polygonable if it admits a polygonation.
\end{defn}

\begin{figure}
\includegraphics[scale=0.6]{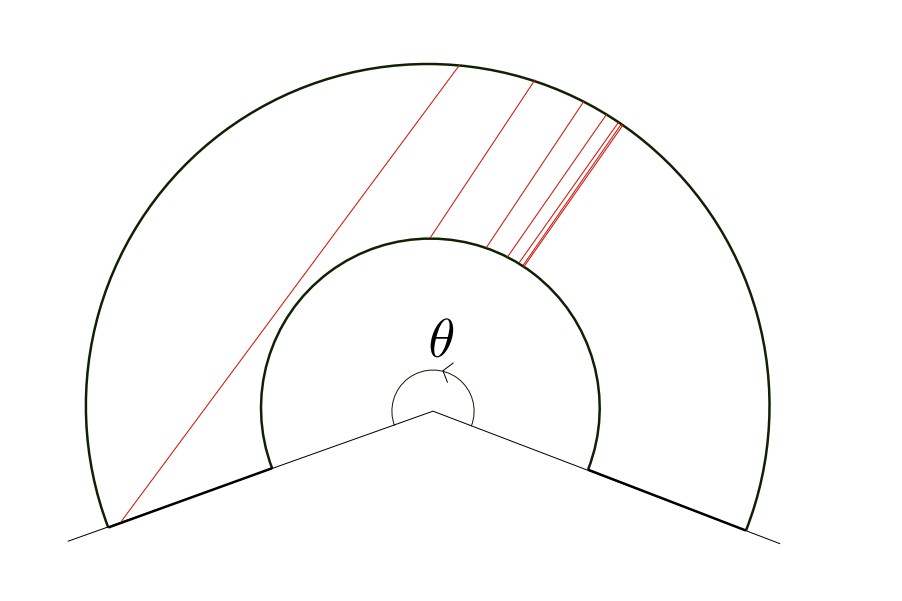}
\caption{A fundamental domain for the action of $z \mapsto 2 z$ on a cone of angle $\theta > \pi$. Any trajectory entering the cylinder is trapped within it forever regardless of the direction of the trajectory, as the one represented here. This property prevents a polygonation to "connect" both sides of the cylinder.
}
\label{figure3}
\end{figure}

Veech's criterion provides a geometric characterization of polygonable surfaces. This theorem is optimal since cylinders of angle at least $\pi$ are not polygonable as shown in Figure~\ref{figure3}. See \cite{DFG19} for a proof.

\begin{thm}\label{thm:Veech}(Veech's criterion)
For a closed dilation surface $\Sigma$ containing at least one singularity, the three following propositions are equivalent:
\begin{itemize}
\item $\Sigma$ is polygonable;
\item $\Sigma$ does not contain a dilation cylinder of angle at least $\pi$;
\item every affine immersion of the open unit disk $\mathbb{D} \subset \mathbb{C}$ in $\Sigma$ extends continuously to its closure $\bar{\mathbb{D}}$.   
\end{itemize}
\end{thm}

\begin{rem}
Up to adding enough singularities of angle $2\pi$ and trivial dilation multiplier, we can nevertheless decompose cylinders of angle at least $\pi$ into smaller cylinders and then into polygons.
\end{rem}

For our purpose, Theorem~\ref{thm:Veech} proves in particular that every dilation surface that is not polygonable carries cylinders and one can focus on polygonable surfaces. 

\subsection{Action of $SL(2,\mathbb{R})$} We now define a natural action of $SL(2,\mathbb{R})$ on the space of dilation surfaces.

\vspace{2mm} \noindent  Let $\Sigma$ be a dilation surface and consider $A \in SL(2,\mathbb{R})$. Let $(U_{i}, \varphi_i)_{i \in I}$ a maximal atlas defining the dilation structure of $\Sigma$. Define $ A\Sigma $ to be the dilation structure defined by the maximal atlas $(U_i, A \circ \varphi_i)_{i \in I}$ where $A$ acts on $\mathbb{C}$ via the standard identification $\mathbb{C} \simeq \mathbb{R}^{2}$. This new atlas indeed defines a dilation structure as $SL(2,\mathbb{R})$ centralizes the group formed by maps $z \mapsto az+b$ where $a \in \mathbb{R}_{+}^{\ast}$ and $b \in \mathbb{C}$.

\vspace{2mm} If the dilation surface was given by gluing a bunch of polygons together, the new surface is also polygonable. Indeed, the image of the initial set of polygons with edges identifications is mapped by the linear action of the matrix $A$ to another set of polygons. Since $A$ is linear, the sides of the polygon that were parallel are still parallel after applying the matrix $A$ so that one can still glue them using dilations of the plane: the resulting dilation surface is the image of $\Sigma$ under the matrix $A$.\newline

A remarkable subgroup of $SL(2,\mathbb{R})$ is formed by matrices $g_{t}=\begin{pmatrix} e^{t} & 0 \\ 0 & e^{-t} \end{pmatrix}$ for $t \in \mathbb{R}$. The flow expands the horizontal direction and contracts the vertical direction.\newline

\section{Delaunay polygonations}\label{sec:Delaunay}

\subsection{Delaunay polygonation}

The goal of this subsection is to define the Delaunay polygonation of a (polygonable) dilation surface. The construction we will give actually is Veech's proof of Theorem~\ref{thm:Veech}. To show that surfaces that do not carry cylinders of angle larger than $\pi$ are polygonable, he proved that the following construction defines a polygonation. We refer to \cite{DFG19} for the full proof. We will here only describe it. \\

The vertices of this polygonation are by definition the singularities of $\Sigma$. The edges of the polygonation are saddle connections: a given saddle connection between singularities $s_{1}$ and $s_{2}$ belongs to the edges of the Delaunay triangulation if there is a closed disk immersed in $\Sigma$ such that $s_{1}$ and $s_{2}$ belong to the boundary circle of this disk and such that there are no other singularities in its interior. 
A disk in $\Sigma$ is said to be Delaunay if it does not contain any singularities in its interior but at least three on its boundary. \\

\begin{figure}
\includegraphics[scale=0.7]{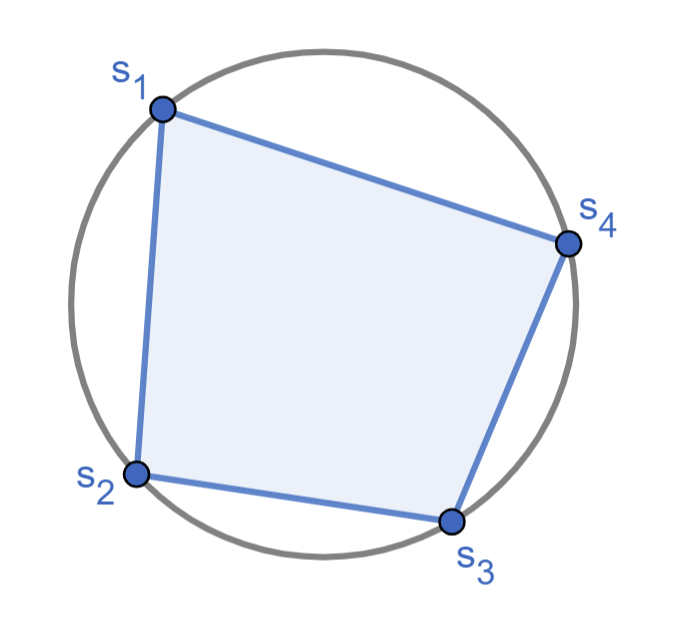}
\caption{A Delaunay disk with four boundary singularities. Their convex hull in the disk is a face of the Delaunay polygonation.}
\label{figure4}
\end{figure}

The faces correspond to what is left after suppressing the edges and the vertices: they are convex polygons whose extremal points belongs all to the same Delaunay disks. Figure~\ref{figure4} illustrates the construction: the disk is a Delaunay disk whose boundary contains singularities $s_{1},s_{2},s_{3},s_{4}$. The quadrilateral is one of the face of the polygonation while its four sides are edges.\\

Note that the Delaunay polygonation gives you a way to recover from an "abstract" polygonable dilation surface a concrete set of polygons that defines it.

\begin{rem}
In this paper, even if surfaces with boundary may appear, we will only consider Delaunay polygonations of dilation surfaces without boundary.
\end{rem}

\subsection{Polygons up to dilation and their limits.}
\label{sec:polygons}
In this paragraph we consider the space of polygons with exactly $p \geq 3$ vertices arising from Delaunay polygonation. We consider these polygons \textbf{as marked and up to dilation}, which means that

\begin{itemize}
\item we think of Delaunay polygons as within the unit circle as we can use a dilation to map the Delaunay circle to the unit one;

\item we keep track of the role of each side and each vertex which is what we mean by marked;

\item a polygon and its image under a rotation are considered to be different (because polygons are considered up to dilation and not similarity).
\end{itemize}

We denote the set described above by $\mathcal{P}_{p}$. Each polygon is characterized by a $p$-uplet $(\theta^{1},\dots,\theta^{p}) \in (\mathbb{R}/2\pi\mathbb{Z})^{p}$ where $\theta^{i}$ is the angle of vertex $i$ in the Delaunay circle and $\theta^{i} \neq \theta^{j}$ for $i \neq j$.

\begin{figure}
\includegraphics[scale=0.8]{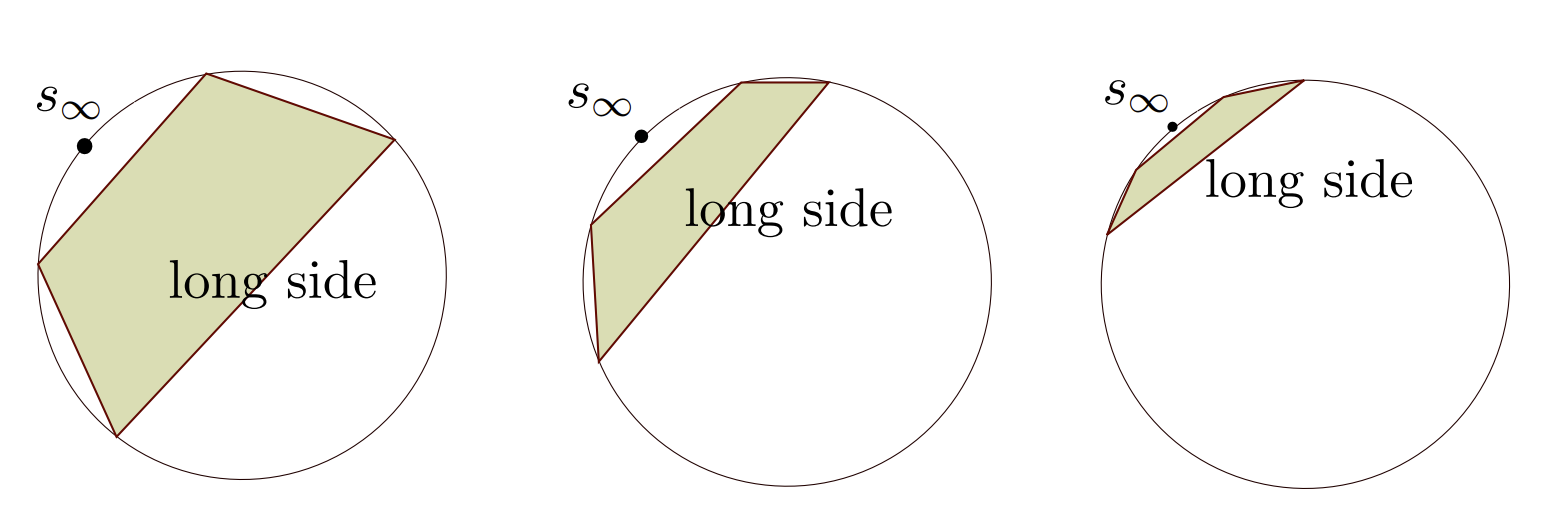}
\caption{The three first polygons of a degenerating sequence of type 1 whose vertices all converge toward $s_{\infty}$.}
\label{figure5}
\end{figure}

\begin{figure}
\includegraphics[scale=0.8]{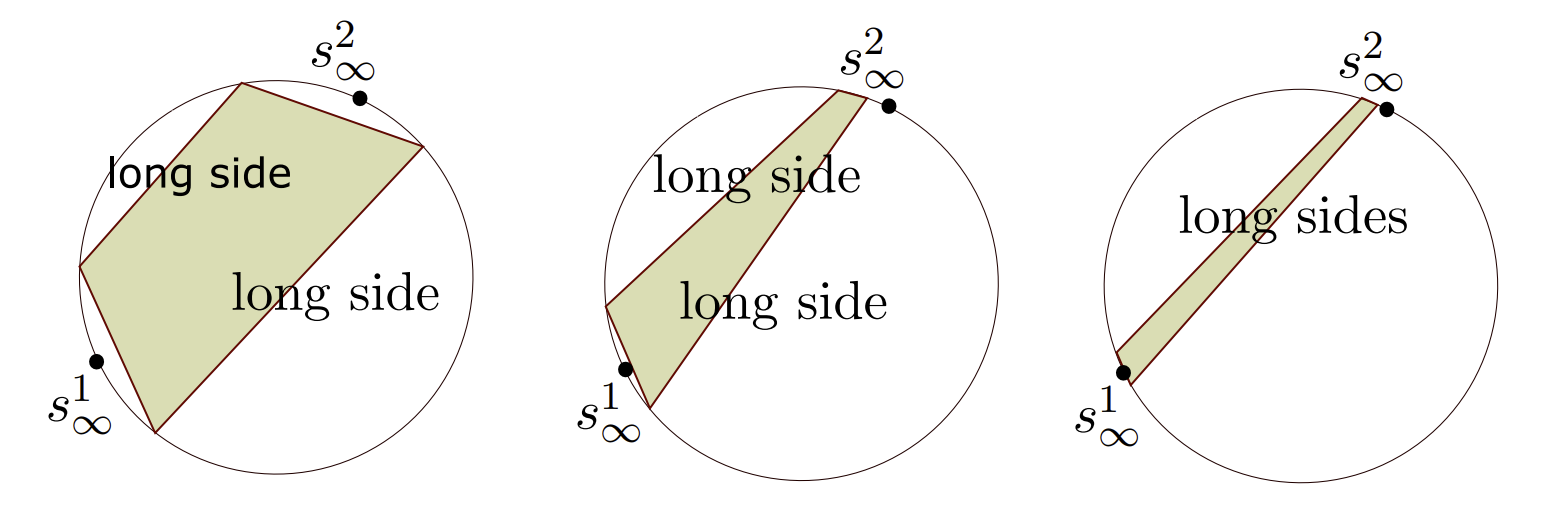}
\caption{The three first polygons of a degenerating sequence of type 2 whose vertices all converge toward either $s^1_{\infty}$ or $s^2_{\infty}$.}
\label{figure6}
\end{figure}

\begin{defn}\label{defn:DCpolygons} Let $p$ be fixed. Consider a sequence of polygons $(P_{n})_{n \in \mathbb{N}}$ in $\mathcal{P}_{p}$. We say that this sequence is \textbf{Delaunay-convergent} if the following conditions hold:
\begin{itemize}
    \item the cyclic ordering of the vertices in the circle is constant;
    \item each vertex $(\theta_{n}^{i})_{n \in \mathbb{N}}$ converges in the circle.
\end{itemize}
Besides, a Delaunay-convergent polygon is:
\begin{itemize}
    \item a \textbf{polygon of type 1} if all the vertices of $(P_{n})_{n \in \mathbb{N}}$ converge towards a given point $s_{\infty}$ of the circle containing all vertices of $P_{n}$ (see Figure~\ref{figure5});
    \item a \textbf{polygon of type 2} if the set of vertices of $(P_{n})_{n \in \mathbb{N}}$ converges towards a set of exactly two points $s_{\infty}^{1}$ and $s_{\infty}^{2}$ of the Delaunay circle (see Figure~\ref{figure6}). The slope of the limit edge in $\mathbb{RP}^{1}$, relating the two remaining vertices, is called the \textbf{limit slope}.
    \item a \textbf{polygon of type 3} if the set of vertices of $(P_{n})_{n \in \mathbb{N}}$ converges towards a set of at least three vertices.
\end{itemize}
\end{defn}

By compactness, one can from any sequence of polygons $(P_{n})_{n \in \mathbb{N}}$ extract a Delaunay-convergent subsequence. \\

We now introduce the following terminology which will be useful when proving our main theorem. 

\begin{defn}\label{defn:longshort}
If $(P_{n})_{n \in \mathbb{N}}$ is of type 1, the longest side of the polygon, corresponding to the closest one from the center of the circle in which it is inscribed, is called the \textbf{long side} while the other sides will be called \textbf{short}. The \textbf{direction} of the polygon in $\mathbb{RP}^{1}$ is the direction given by its long side \newline
If $(P_{n})_{n \in \mathbb{N}}$ is of type 2 or 3, the longest sides whose vertices converge to different limit point will be called the \textbf{long sides} while the others sides are called \textbf{short}. In the type 2 case, the directions tangent to the circle at the two remaining vertices are called the \textbf{short side limit slopes} as the short sides are asymptotic to this direction.
 \end{defn}

Be aware that the terminology is about sequences of polygons and more precisely about their asymptotic behaviour: one can change finitely many polygons of the sequence without changing its long or short sides.  \\
 
By an harmless abuse of notation, we will refer to a sequence of polygons $(P_{n})_{n \in \mathbb{N}}$ in a Delaunay-convergent sequence as a \textbf{polygon}. We will use the terms of degenerating polygons. The terms long sides and short sides for these polygons will refer similarly to sequences of edges.

\subsection{Delaunay-convergent sequences of dilation surfaces}\label{sec:sequence}

In this paragraph we consider sequences of dilation surfaces of fixed topological type (the underlying marked topological surfaces are isomorphic).\newline

Let $(\Sigma_{n})_{n \in \mathbb{N}}$ be a sequence of dilation surfaces of same topological type. Up to extracting a subsequence, we can assume that their Delaunay polygonations are all combinatorially equivalent. Precisely this means that for any $n \in \mathbb{N}$ their Delaunay polygonations have the same pattern. We label for each $n$ the set $I$ of polygons $(P_{i,n})_{i \in I}$, in such a way that 
\begin{itemize}
 \item the sequence $(P_{i,n})_{n \in \mathbb{N}}$ has always the same numbers of sides;
 \item one can mark the sides of the polygons in a way that the gluing pattern of the sides of the marked polygons $(P_{i,n})_{n \in \mathbb{N}}$ is constant with respect to the marking.
\end{itemize}
In that case we say that sequence of surfaces $(\Sigma_{n})_{n \in \mathbb{N}}$ has \textbf{constant Delaunay pattern}.

\begin{defn}\label{defn:DCdilation}
A sequence $(\Sigma_{n})_{n \in \mathbb{N}}$ is said to be \textbf{Delaunay-convergent} if:
\begin{enumerate}
    \item sequence $(\Sigma_{n})_{n \in \mathbb{N}}$ has constant Delaunay pattern;
    \item every polygon $(P_{i,n})_{n \in \mathbb{N}}$ is Delaunay-convergent (see Definition~\ref{defn:DCpolygons}).
\end{enumerate}
We refer to the edges of these polygons as the \textbf{Delaunay edges} of the pattern.
\end{defn}

For a given sequence of dilation surfaces of fixed topological type, there are finitely many Delaunay patterns so we can always extract a Delaunay-convergent subsequence.

\subsection{Maximal domains of type 1}\label{sec:domain1}

Properties of Delaunay polygonations induce constraints on the Delaunay patterns involving polygons of type 1.

\begin{prop}\label{prop:longside}
In a Delaunay-convergent sequence $(\Sigma_{n})_{n \in \mathbb{N}}$ of polygonable dilation surfaces, the long side $(L_{n})_{n \in \mathbb{N}}$ of a polygon $(P_{n})_{n\in \mathbb{N}}$ of type 1 can only be incident to a short side of a polygon (of any type).
\end{prop}

\begin{proof}
We assume that in addition to be a side of $(P_{n})_{n\in \mathbb{N}}$, edge $(L_{n})_{n\in\mathbb{N}}$ is a long side of a polygon $(Q_{n})_{n\in \mathbb{N}}$. In these cases, vertices of $(P_{n})_{n\in \mathbb{N}}$ distinct from the ends of $(L_{n})_{n\in\mathbb{N}}$ are included for $n$ large enough in the interior of the Delaunay disk of polygons $(Q_{n})_{n\in \mathbb{N}}$, see Figure~\ref{figure7}. This is a contradiction.
\begin{figure}[h!]
\includegraphics[scale=0.6]{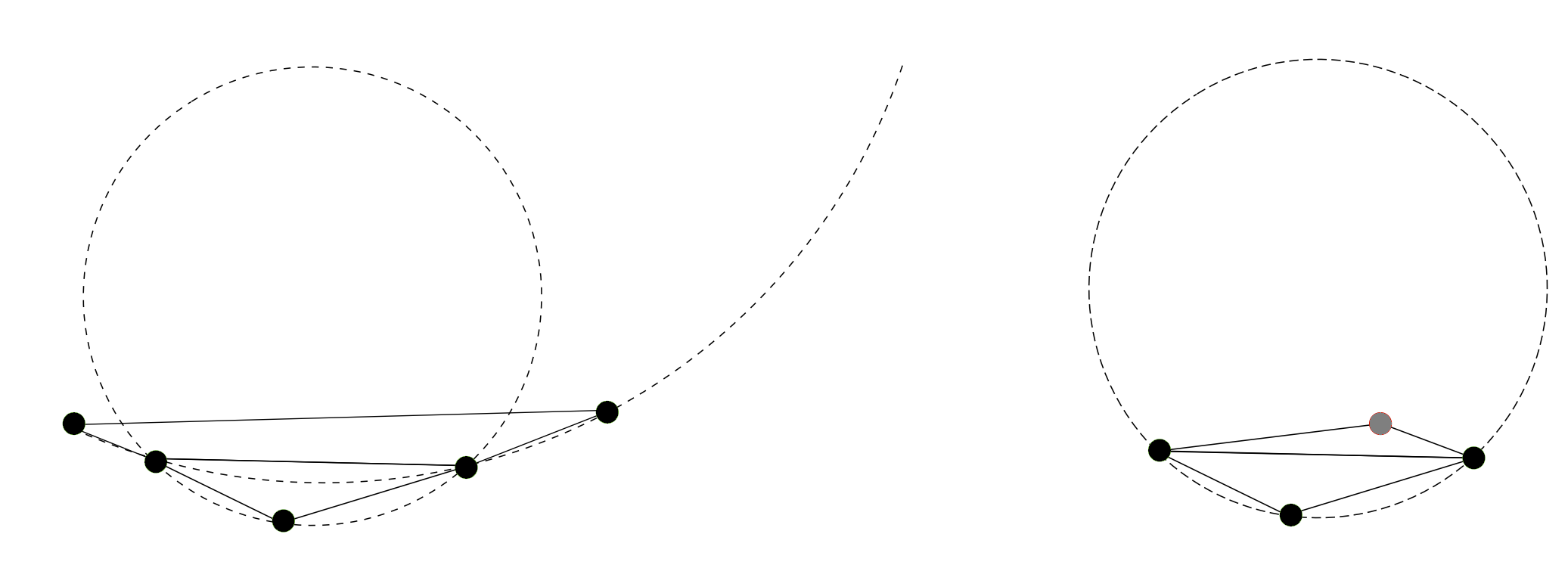}
\caption{On the left, a configuration of two polygons of type 1: none of the singularities lies inside the two Delaunay disks. As the polygons shrink, the Delaunay disk tend to cover half a plane. On the right a forbidden configuration: the grey dot lies inside a large Delaunay disk.}
\label{figure7}
\end{figure}
\end{proof}

Following Proposition~\ref{prop:longside}, the long side of a polygon of type 1 is always incident to the short side of another polygon. We say that two polygons of type 1 belong to the same \textbf{domain of type 1} if the long side of the first one is incident to a short side of the second one. The equivalence relation generated by these relations define classes. This way, each polygon of type 1 belongs to a unique \textbf{maximal domain of type 1}.
\par
Besides, in a maximal domain of type 1, any internal edge is a long side a polygon while being a short side of another (it may happen that the polygons coincide). Thus, the incidence graph of a maximal domain is actually an oriented graph with a unique oriented edge leaving each vertex (since each polygon of type 1 has a unique long side). It follows from that there are two types of maximal domains of type 1:
\begin{itemize}
    \item \textbf{non-cyclic domains} where the incidence graph is a rooted tree (the edges being oriented towards the root);
    \item \textbf{cyclic domains} where the oriented incidence graph contains a unique (oriented) cycle.
\end{itemize}
Since a maximal domain of type 1 is connected, there is no other type of graphs of incidence.

\subsection{Maximal domains of type 2}\label{sec:chain2}

Polygons of type 2 have two long sides. Two polygons of type 2 glued along an edge that is a long side for each of them belong to a same domain of type 2. This way, each polygon of type 2 belongs to a unique \textbf{maximal domain of type 2}. These domains can be \textbf{cyclic} or not. In the non-cyclic case we call the two long edges that are not glued with another polygon of type 2 the \textbf{extremal edges}. \newline

\begin{rem}
The case of a long side of a polygon of type 2 glued along a short side of another polygon of type 2 can happen. We can obtain such a configuration in a variant of the degeneration presented in Figure~\ref{figure9}. If the modulus of the connecting flat cylinder decreases to zero (instead of going to infinity), the cylinder is a maximal domain of type 2 and its upper extremal edge is glued on the short side of a polygon of type 2.
\end{rem}

\begin{figure}[h]
\includegraphics[scale=0.8]{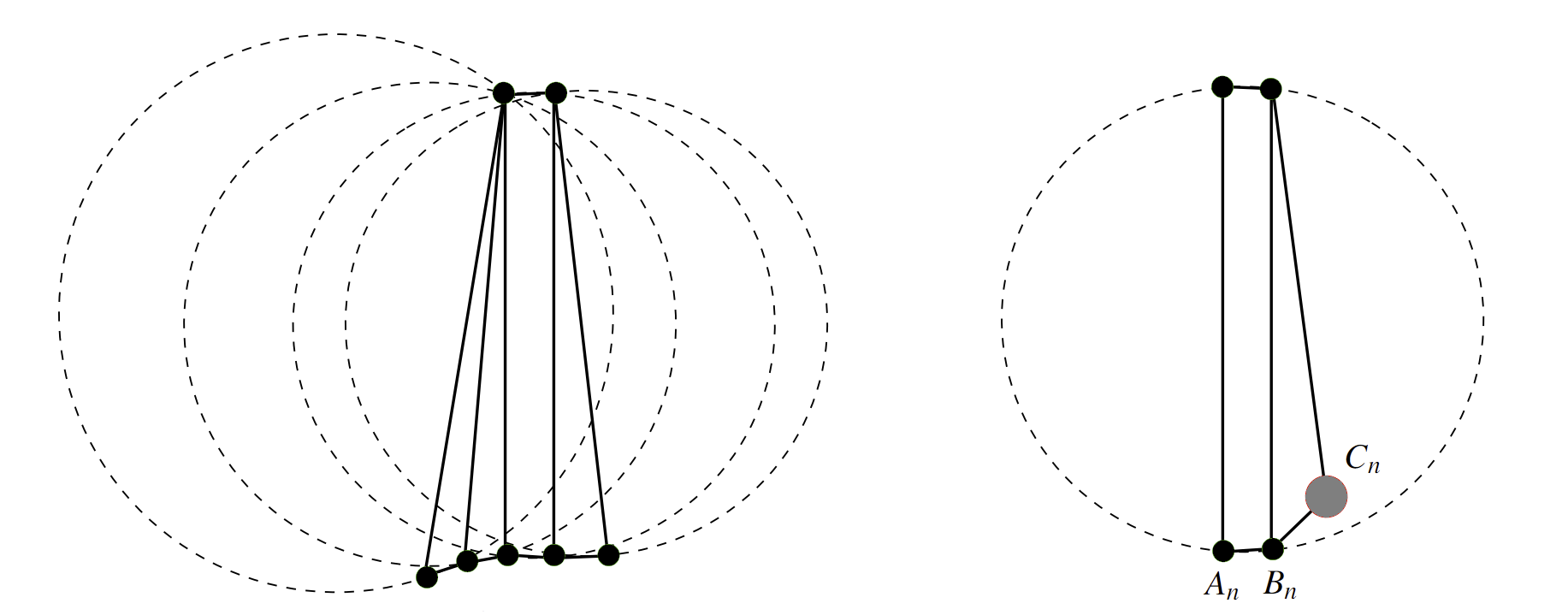}
\caption{On the left, an admissible configuration of four degenerating polygons of type 2: none of the singularities lies inside the Delaunay disks. As the polygons shrink, the Delaunay disks tend to cover a half-plane. On the right a forbidden configuration: the grey dot lies inside a large Delaunay disk.}
\label{figure8}
\end{figure}

An observation that will be needed in the proof of Theorem~\ref{thm:main}, is that short sides of maximal domains of type 2 form "concavely shaped" curves, as shown in Figure~\ref{figure8}. It proceeds from the following statement.

\begin{prop}\label{prop:concave}
We consider a Delaunay-convergent sequence $(\Sigma_{n})_{n \in \mathbb{N}}$ of polygonable dilation surfaces. In the polygon formed by two degenerating polygons of type 2 or 3 glued along a common long side, the magnitude of the limit inner angle between two consecutive short sides is at least $\pi$.
\end{prop}

\begin{proof}
Two consecutive short sides belonging to the same polygon of type 2 or 3 have the same limit slope because their endpoints converge to the same limit point in the Delaunay circle (see Definition~\ref{defn:longshort}). Therefore, the limit inner angle between them is equal to $\pi$.
\par
Now we consider the case of two consecutive short sides $[A_{n},B_{n}]_{n \in \mathbb{N}}$ and $[B_{n},C_{n}]_{n \in \mathbb{N}}$ belonging respectively to two distinct incident degenerating polygons $(P^{1}_{n})_{n\in \mathbb{N}}$ and $(P^{2}_{n})_{n\in \mathbb{N}}$ of type 2 or 3. These two sides have well-defined limit slopes (corresponding to the slope of the tangent line at their limit point in their Delaunay circle). We will assume by contradiction that the limit inner angle $\theta$ between these sides at $B_{n}$ is strictly smaller than $\pi$. \\

It follows that there is $N>0$ such that for any $n \geq N$, the segment $[B_{n}C_{n}]$ intersects the Delaunay disk $\mathcal{D}^{1}_{n}$ of $P^{1}_{n}$. Since by hypothesis $C_{n}$ cannot belong to $P^{1}_{n}$, the segment $[B_{n}C_{n}]$ intersects the boundary of $P^{1}_{n}$ in some point $C_{n}'$. The triangle $A_{n}B_{n}C_{n}'$ is inscribed in the Delaunay circle that bounds $\mathcal{D}^{1}_{n}$. \\

Since $P_{n}^{1} \cup P_{n}^{2}$ is contractible, it can be endowed with a flat metric in such a way that $[A_{n},B_{n}]_{n \in \mathbb{N}}$ and $[B_{n},C_{n}]_{n \in \mathbb{N}}$ have meaningful lengths. The latter metric is normalized by fixing the radius of Delaunay disk $\mathcal{D}^{1}_{n}$ to $1$. As $n \rightarrow + \infty$, the length of $[A_{n},B_{n}]$ shrinks to zero. Since the inner angle at $B_{n}$ converges to $\theta<\pi$, the length of $[B_{n}C'_{n}]$ converges to some nonzero limit. It follows that the length of $[B_{n}C_{n}]$ cannot decrease to zero as $n$ tends to infinity. In 
 $(P^{2}_{n})_{n\in \mathbb{N}}$, the length of $[B_{n}C_{n}]$ does not become negligible in comparison with the length of the common edge between $(P^{1}_{n})_{n\in \mathbb{N}}$ and $(P^{2}_{n})_{n\in \mathbb{N}}$. In other words, $[B_{n}C_{n}]$ is not a short side of $(P^{2}_{n})_{n\in \mathbb{N}}$ and we get a contradiction.
\end{proof}

\subsection{Delaunay limits}\label{sec:DL}

\begin{figure}[h!]
\includegraphics[scale=0.8]{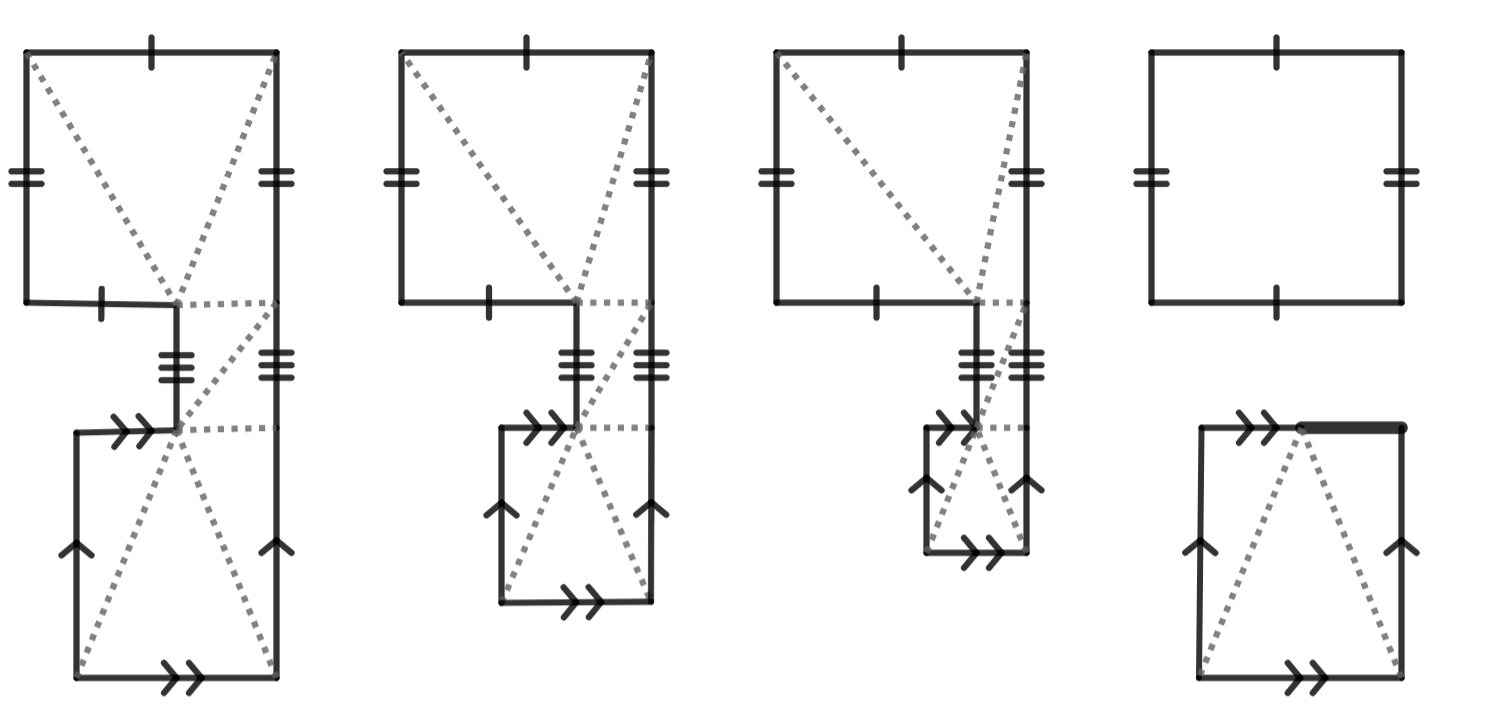}
\caption{The three first drawings represent some terms of a degenerating dilation surface of genus two with one singularity whose 'connecting' flat cylinder degenerates as its modulus goes to infinity.
At the level of the Delaunay polygonation, the two parts left after removing the connecting cylinder converge to:\newline
(i) a flat torus which is the "Delaunay limit" of the upper part;\newline
(ii) a torus with one horizontal boundary component (in bold). 
}
\label{figure9}
\end{figure}

For any Delaunay-convergent sequence $(\Sigma_{n})_{n \in \mathbb{N}}$ of closed dilation surfaces, we can define a \textbf{Delaunay limit} $\Sigma_{\infty}$ formed by the polygons that do not completely degenerate, see Figure~\ref{figure9} for an example.\newline

Limit surface $\Sigma_{\infty}$ will be a polygonable dilation surface. However, we should be careful. The limit surface can have several connected components. It can also have a boundary and it can even be empty.

\begin{defn}\label{defn:DL}
Let $(\Sigma_{n})_{n \in \mathbb{N}}$ be a Delaunay-convergent sequence of closed dilation surfaces. We define the \textbf{Delaunay limit} $\Sigma_{\infty}$ in the following way:
\begin{itemize}
    \item $\Sigma_{\infty}$ is the union of limits of polygons of type 3 in $(\Sigma_{n})_{n \in \mathbb{N}}$;
    \item two sides of limit polygons of $(\Sigma_{n})_{n \in \mathbb{N}}$ are identified to each other if they are incident in the Delaunay pattern of the sequence;
    \item the sides of limit polygons are also identified to each other if they are connected by a non-cyclic maximal domain of type 2 (see Section~\ref{sec:chain2}).
\end{itemize}
As they degenerate to polygons with empty interior, polygons of type 1 and 2 completely disappear in the Delaunay limit $\Sigma_{\infty}$ (see Figures~\ref{figure9} and~\ref{figure10}).
\end{defn}

The simple, but key, feature about this notion of limit is that 'carrying a cylinder' is an open property. Namely, if a sequence of dilation surfaces has a Delaunay limit which carries a cylinder then for any large enough index of the sequence the corresponding surface also carries a cylinder. 

\begin{prop}\label{prop:cylinderOPEN}
We consider a Delaunay-convergent sequence $(\Sigma_{n})_{n \in \mathbb{N}}$. If its Delaunay limit $\Sigma_{\infty}$ is non-empty and contains a closed geodesic in some direction $\theta \in \mathbb{RP}^{1}$, then for any $\epsilon>0$, there is $N>0$ such that for any $n \geq N$, $\Sigma_{n}$ contains a closed geodesic in a direction of $]\theta-\epsilon,\theta+\epsilon[$.
\end{prop}

\begin{proof}
We denote by $\gamma$ a closed geodesic of slope $\theta$ in the limit surface $\Sigma_{\infty}$. Such a geodesic must cross an edge of the Delaunay polygonation as there is no closed geodesic contained in the interior of a Delaunay polygon. Let us denote by $[A,B]$ an edge crossed by $\gamma$ and by $[A,B]_\gamma$ the intersection of $[A,B]$ with $\gamma$. By definition, closed geodesics do not contain any singularity. It follows that $[A,B]_{\gamma}$ belongs to the interior of $[A,B]$ (it is not a singularity). Moreover, as $[A,B]_{\gamma}$ belongs to a periodic leaf of $\mathcal{F}_{\theta}$, the foliation $\mathcal{F}_{\theta}$ on $\Sigma_{\infty}$ has a well-defined first return map on a neighbourhood $[C,D]$ of $[A,B]_{\gamma}$, see Figure~\ref{figure12}. \\

By definition of the Delaunay limit, the edge $[A,B]$ is the limit edge of a sequence of long sides $([A_n,B_n])_{n \in \mathbb{N}}$ of a polygon of type 3 in the Delaunay polygonation of $(\Sigma_n)_{n \in \mathbb{N}}$. As the polygons converge, the unique (up to translation) complex affine mapping $\mathcal{A}_n$ of the plane that maps $[A_n,B_n]$ to $[A,B]$ converges to the identity as $n \to \infty$. We set $x_n := \mathcal{A}_n^{-1}([A,B]_{\gamma})$. Note that $x_n \to [A,B]_{\gamma}$ as $n \to \infty$. 
\par
\begin{figure}[h!]
\includegraphics[scale=0.4]{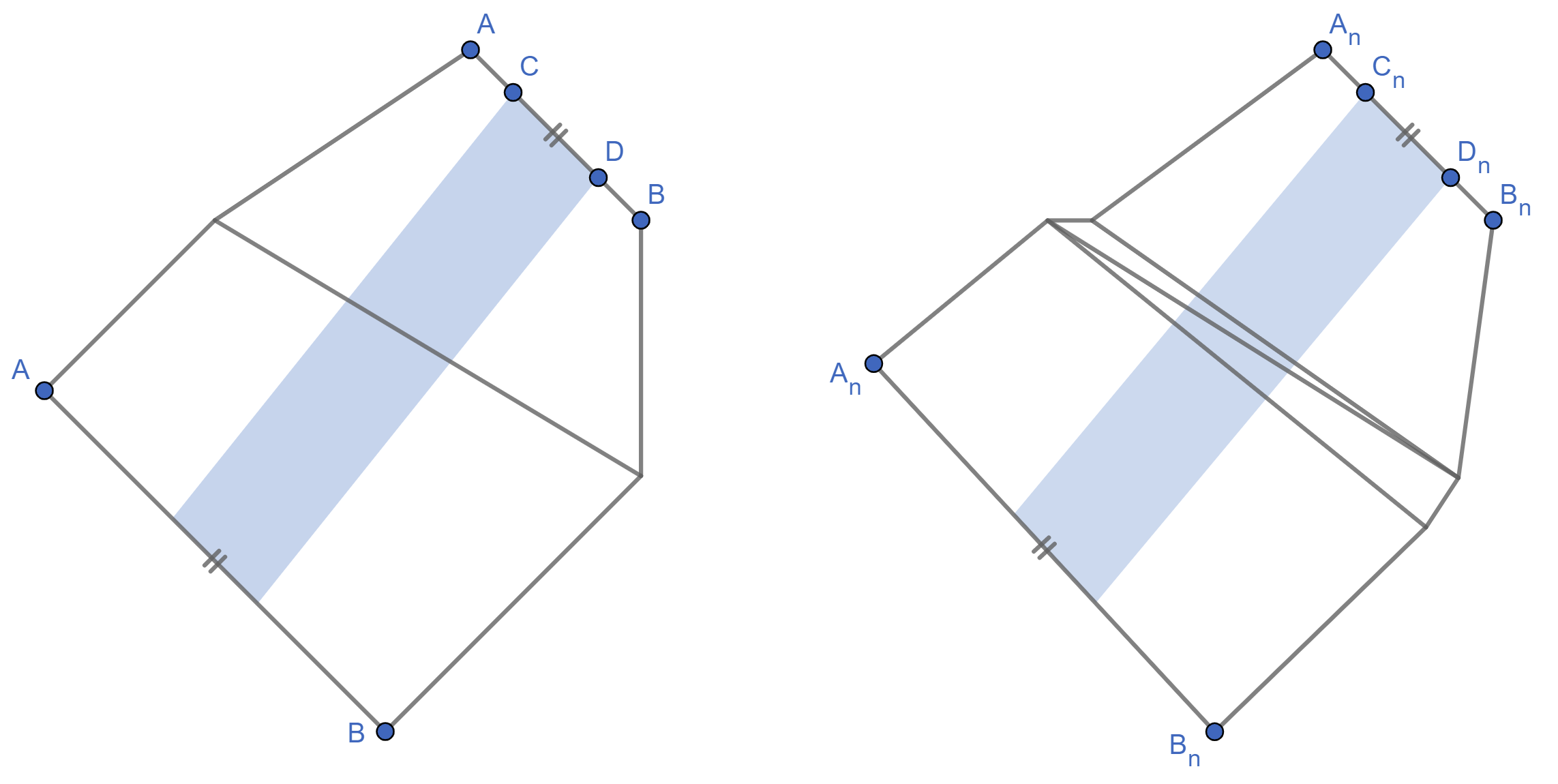}
\caption{The bands here are represented in shadowed blue. On the left it is the neighbourhood of some periodic orbit and on the right the corresponding band, defined by using the affine identification $\mathcal{A}_n$ of the edges.}
\label{figure12}
\end{figure}

For $n$ large enough, the leaf of the foliation $\mathcal{F}_{\theta}$ in $\Sigma_{n}$ starting at $x_n$ crosses the edges corresponding to the edges crossed by $\gamma$ in $\Sigma_{\infty}$ (several edges of $\Sigma_{n}$ can correspond to the same edge of $\Sigma_{\infty}$ if they are long sides of the same non-cyclic maximal domain of type 2) and then crosses back $[A_n,B_n]$ at some point. \\

As the finitely many polygons encountered converge toward non-degenerate polygons or degenerates toward 'edges', there is a bound $N>0$ such that for any $n \geq N$, the first return map $T_{n}$ of $[A_n,B_n]$ is well defined on a neighbourhood $[C_n, D_n] := \mathcal{A}_n^{-1}([C,D])$. \\

All the oriented leaves of $\mathcal{F}_{\theta}$ starting from $[C_n, D_n]$, taken up to their first return on $[A_n,B_n]$, give rise to a band $\mathcal{B}_{n}$ (a parallelogram) whose sides contained in $[A_{n},B_{n}]$ partially coincide. In particular, $\mathcal{B}_{n}$ contains a closed geodesic. Similarly, we define $\mathcal{B}_{\infty}$ in $\Sigma_{\infty}$. \\

For an arbitrarily small $\eta>0$, we can choose a neighbourhood $[C,D]$ of $[A,B]_{\gamma}$ such that the slopes of the two diagonals of $\mathcal{B}_{\infty}$ are contained in $]\theta-\eta,\theta+\eta[$. Then, provided $n$ is large enough, the slopes of the diagonals of $\mathcal{B}_{n}$ can be made arbitrarily close to slopes of the diagonals of $\mathcal{B}_{\infty}$. The slope of a closed geodesic contained in $\mathcal{B}_{n}$ belongs to an interval whose ends are the slopes of the diagonals of $\mathcal{B}_{n}$
It follows that for any $\epsilon>0$, there is $N>0$ such that for any $n \geq N$, $\mathcal{B}_{n}$ contains a closed geodesic whose slope is contained in $]\theta-\epsilon,\theta+\epsilon[$.
\end{proof}

We also need to keep track of the polygons involved in the Delaunay limit. To this purpose, we introduce the notion of \textbf{core sequence}.

\begin{defn}\label{defn:core}
For a given Delaunay-convergent sequence $(\Sigma_{n})_{n \in \mathbb{N}}$, the \textbf{core sequence} $(\mathcal{C}\Sigma_{n})_{n \in \mathbb{N}}$ is defined for each $n$ as the union of:
\begin{itemize}
    \item polygons of type 3;
    \item non-cyclic maximal domains of type 2 in which at least one extremal edge is incident to a long side of a polygon of type 3.
\end{itemize}
\end{defn}

\begin{rem}\label{rem:CC}
It follows from Definitions~\ref{defn:DL} and~\ref{defn:core} that each connected component of $\Sigma_{\infty}$ corresponds to a unique connected component of $\mathcal{C}\Sigma_{n}$. Boundary saddle connections of $\Sigma_{\infty}$ correspond to long boundary edges of the core.
\end{rem}

\begin{figure}[h!]
\includegraphics[scale=0.6]{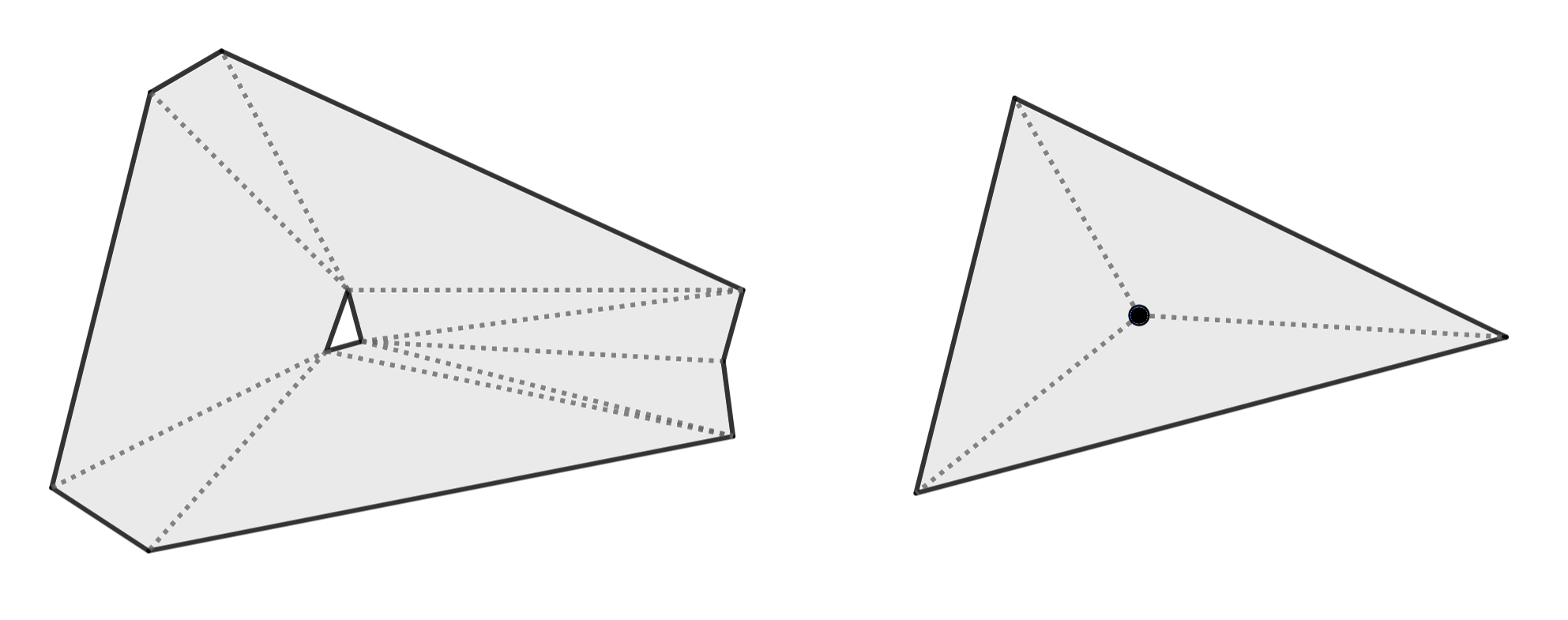}
\caption{On the left, surface $\mathcal{C}\Sigma_{n}$ of the core sequence with two boundary components. The exterior boundary component contains three long sides while the interior boundary component is formed by short sides only. On the right, Delaunay limit $\Sigma_{\infty}$ with an exterior boundary formed by three sides and a singularity at the center.}
\label{figure10}
\end{figure}

Maximal domains of type 1, maximal domains of type 2 and connected components of the core are the fundamental pieces of the decomposition we will use in the proof of Theorem~\ref{thm:main}.

\begin{defn}\label{defn:PIECES}
Polygons of $(\Sigma_{n})_{n \in \mathbb{N}}$ are grouped into \textbf{Delaunay pieces} that are:
\begin{itemize}
    \item the connected components of the core $(\mathcal{C}\Sigma_{n})_{n \in \mathbb{N}}$;
    \item the maximal domains of type 1;
    \item the maximal domains of type 2 that do not belong to the core.
\end{itemize}
It follows from Definition~\ref{defn:core} that every polygon of $(\Sigma_{n})_{n \in \mathbb{N}}$ belongs to a unique Delaunay piece.
\end{defn}

The important feature of the decomposition into Delaunay pieces is that their boundary edges cannot be long sides for both of their incident polygons.

\begin{lem}\label{lem:Delpiece}
If a Delaunay edge $(L_{n})_{n \in \mathbb{N}}$ of $(\Sigma_{n})_{n \in \mathbb{N}}$ is a long side for its two incident polygons, then it cannot belong to the boundary of a Delaunay piece.
\end{lem}

\begin{proof}
It follows from Lemma~\ref{prop:longside} that none of the incident polygons of $(L_{n})_{n \in \mathbb{N}}$ can polygons of type 1. If both of them are polygons of type 2, then they belong to the same maximal domain and $(L_{n})_{n \in \mathbb{N}}$ is not a boundary edge of a Delaunay piece. If at least one of the two polygons incident to $(L_{n})_{n \in \mathbb{N}}$ is of type 3, then it follows from Definition~\ref{defn:core} that $(L_{n})_{n \in \mathbb{N}}$ is an interior edge of some connected component of the core $(\mathcal{C}\Sigma_{n})_{n \in \mathbb{N}}$.
\end{proof}

\section{Overview of proof of Theorem~\ref{thm:main}}\label{sec:overview}

As mentioned above, we will argue by induction on the complexity of the closed surface $\Sigma$ (see Definition~\ref{defn:complexity}). It is easy to deal with the case of smallest complexity, the case of flat tori formed by a pair of triangles. The induction assumption is that any closed dilation surface of complexity lower than $k$ carries cylinders in a dense set of directions. We want to show that it is still the case for surfaces of complexity $k+1$.\\

It follows from Lemma~\ref{lem:cone} that for a dilation surface $\Sigma$ without boundary, an open set of $\mathbb{RP}^{1}$ that does not contain any direction of saddle connection contains the direction of a closed geodesic. Therefore, it remains to prove that any direction $\theta \in \mathbb{RP}^{1}$ that is approached by directions of saddle connections is also approached by directions of closed geodesics. In other words, any open subset of $\mathbb{RP}^{1}$ containing the direction of a saddle connection should also contain the direction of a closed geodesic. Up to the action of an element of $\mathrm{SL}_2(\RR)$, we can assume that $\Sigma$ contains a vertical saddle connection $\gamma$ and $U$ is an open subset of $\mathbb{RP}^{1}$ containing the vertical direction.\\ 

To rely on the induction assumption, we will use the Teichmüller flow $g_t$ that contracts the vertical direction, and therefore the saddle connection $\gamma$, and that expands the horizontal direction. We have seen in Section~\ref{sec:sequence} that one can extract a sequence of times $t_n \to + \infty$ such that the sequence $(\Sigma_{n})_{n \in \mathbb{N}} := (g_{t_n}\Sigma)_{n \in \mathbb{N}}$ is Delaunay-convergent. We are doing induction on closed surfaces and this is why the case where $\Sigma_{\infty}$ has a boundary will need to be handled separately. We first rule out the case in which $(\Sigma_{n})_{n \in \mathbb{N}}$ Delaunay-converges toward a closed surface $\Sigma_{\infty}$ of same complexity (see Section~\ref{sec:DL} for a precise definition of Delaunay limits). 

\begin{prop}
\label{prop.degenerating}
Let $\Sigma$ be a closed dilation surface that carries a vertical saddle connection and a sequence of times $t_n \to +\infty$ such that $(g_{t_n}\Sigma)_{n \in \mathbb{N}}$ Delaunay-converges toward a surface $\Sigma_{\infty}$. Then one of the following statements holds:
    \begin{itemize}
        \item for any $\epsilon > 0$, $\Sigma$ carries a  cylinder whose direction belongs to $] \frac{\pi}{2} - \epsilon, \frac{\pi}{2} + \epsilon[$;
        \item $\Sigma_\infty$ is of strictly smaller complexity than $\Sigma$. 
    \end{itemize} 
\end{prop}

Proposition~\ref{prop.degenerating} is proved in Section~\ref{sec:divergent}. Note that in the first case there is nothing more to be proven. In the second case, as soon as $\Sigma_\infty$ is non-empty and contains a component without boundary, we can also conclude. Indeed, the induction assumption guarantees that directions of closed geodesics of $\Sigma_{\infty}$ are dense in $\mathbb{RP}^{1}$. In particular, $\Sigma_{\infty}$ contains closed geodesics whose direction are arbitrarily close to the vertical direction.
\par
One can now conclude using Proposition~\ref{prop:cylinderOPEN} that provided $n$ is large enough, $\Sigma_{n}$ contains a closed geodesic whose direction is arbitrarily close to the the vertical direction. It follows that $\Sigma=g_{t_{n}}\Sigma_{n}$ contains a closed geodesic whose direction is even closer to the vertical direction.
\par
It then remains to deal with two cases:
\begin{itemize}
    \item $\Sigma_{\infty}$ is empty;
    \item every connected component of $\Sigma_{\infty}$ has a nonempty boundary.
\end{itemize}
We will deal with these two cases at once by thoroughly examining how the Delaunay polygonations, and especially the Delaunay pieces, degenerate under the Teichmüller flow. \\

The key notion here is the one of \textbf{short and long boundary}. We will say that a Delaunay piece (see Definition~\ref{defn:PIECES}) \textbf{has a long boundary side} if one of its boundary sides is the long boundary (in the sense of Definition~\ref{defn:longshort}) of a Delaunay polygon that belongs to the Delaunay piece. \\

Recall that a Delaunay piece is either a component of the core, or a maximal domain of type 1 or 2. A boundary side of the core can be short or long, but the short ones disappear in $\Sigma_{\infty}$ by construction. In particular, if the limit of a Delaunay piece that belongs to the core has a boundary side then the Delaunay piece in question must have a long boundary side. A maximal domain of type 1 or 2 must have, by construction, a long boundary side except in the case where it is cyclic. We summarize the content of this discussion within the following structural lemma.

\begin{lem}\label{lem:DegenerationPiece}
Let $\Sigma$ be a dilation surface and $t_n \to +\infty$ such that $(g_{t_n}\Sigma)_{n \in \mathbb{N}}$ Delaunay-converges toward a surface $\Sigma_{\infty}$ of strictly smaller complexity than $\Sigma$. Then at least one the following conditions holds:
\begin{enumerate}
\item one of the Delaunay pieces converges toward a closed non-empty dilation surface of smaller complexity;
\item one of the Delaunay pieces is a cyclic maximal domain of type 1 or of type 2;
\item all the Delaunay pieces have at least one long boundary side. 
\end{enumerate}
\end{lem}

As mentioned above, the first case is dealt with using the induction assumption. The second case will follow from the following result that will be proven in Section~\ref{sec:control}.

\begin{prop}
\label{prop.maximaldomains}
Let $\Sigma$ be a dilation surface and $t_n \to +\infty$ such that $(g_{t_n}\Sigma)_{n \in \mathbb{N}}$ Delaunay-converges and such that at least one of its Delaunay pieces is a cyclic maximal domain of type 1 or of type 2. Then for any $\epsilon> 0$, $\Sigma$ carries a cylinder whose direction belongs to $]\frac{\pi}{2} - \epsilon, \frac{\pi}{2} + \epsilon[$.
\end{prop}

One is then left to analyse the last case, in which all the Delaunay pieces have at least one long boundary side. This is the most subtle part of the article. 

\begin{prop}
\label{prop.petiteportegrandeporte}
Let $\Sigma$ be a dilation surface and $t_n \to +\infty$ such that $(g_{t_n}\Sigma)_{n \in \mathbb{N}}$ Delaunay-converges and such that all the Delaunay pieces in $(g_{t_{n}}\Sigma)_{n \in \mathbb{N}}$ have at least one long boundary side.
\par
Then, for any open set $U \subset \mathbb{RP}^{1}$, there is $N>0$ such that for any $n \geq N$, $\Sigma_n$ contains closed geodesics whose direction belong to $U$. 
\end{prop}

The above proposition shows in particular that $\Sigma$ carries closed geodesics whose directions are as close as we want to the vertical one as, as usual, non-horizontal closed geodesics of $\Sigma_n$ are images of almost vertical ones of $\Sigma$ under the Teichmüller flow. This proves that any open set of $\mathbb{RP}^{1}$ containing the vertical direction contains the direction of a closed geodesic of $\Sigma$. The three cases of Lemma~\ref{lem:DegenerationPiece} are settled. Theorem~\ref{thm:main} is now proven.\\

The proof of Proposition~\ref{prop.petiteportegrandeporte} will be given in Section~\ref{sec.petiteportegrandeporte}.

\section{The non-degenerating case (Proof of Proposition~\ref{prop.degenerating})}
\label{sec:divergent}

In this Section, we apply the Teichm\"{u}ller flow to a closed dilation surface $\Sigma$ containing a vertical saddle connection. For a sequence of times $t_n \to + \infty$ such that $(g_{t_n}\Sigma)_{n \in \mathbb{N}}$ Delaunay-converges to a limit surface $\Sigma_{\infty}$ with the same complexity as $\Sigma$, we prove that $\Sigma_{\infty}$ (and subsequently $\Sigma$) contains closed geodesics whose directions are arbitrarily close to $\frac{\pi}{2}$.\\

\begin{rem}
Before entering the proof of the above proposition, let us mention that there is a class of dilation surfaces, called quasi-Hopf surfaces (see \cite{Ta21} for details), having a vertical saddle connection and whose Teichmüller orbit is periodic. These surfaces decompose into disjoint dilation cylinders whose boundary saddle connections are either horizontal or vertical. We cannot extract from the Teichmüller orbit of these surfaces a sequence that Delaunay-converges to a surface of smaller complexity.
\end{rem}

We recall that for positive times, the Teichm\"{u}ller flow expands the horizontal direction and contracts the vertical direction. We first prove the existence of a vertical saddle connection $\gamma_{\infty}$ in $\Sigma_{\infty}$.

\begin{lem}\label{lem.demdivergedichotomy}
Let $\Sigma$ be a closed dilation surface that carries a vertical saddle connection and a sequence of times $t_n \to +\infty$ such that $(g_{t_n}\Sigma)_{n \in \mathbb{N}}$ Delaunay-converges toward a surface $\Sigma_{\infty}$ having the same complexity as $\Sigma$. Then the limit surface $\Sigma_{\infty}$ carries at least one vertical saddle connection.
\end{lem}

\begin{proof} We first prove that for $n$ large enough the vertical saddle connection of $g_{t_{n}}\Sigma$ belongs to an edge of the Delaunay polygonation. Indeed, in the dilation surface $\Sigma$, there is an affine immersion of an elliptic domain $D$ of eccentricity $e<1$ such that vertical saddle connection $\gamma$ coincides with the image of the major axe. The ratio of lengths between the (horizontal) semi-minor axis and the (vertical) semi-major axis is $\sqrt{1-e^{2}}$. Teichm\"{u}ller flow will deform the ellipse. For $T=-\frac{ln(1-e^{2})}{4}>0$, saddle connection $g_{T}\gamma$ in surface $g_{T}\Sigma$ is the vertical diameter of immersed disk $g_{T}D$. Consequently, for any $t \geq T$, $g_{t}\gamma$ is an edge of the Delaunay polygonation of dilation surface $g_{t}\Sigma$. \\

We now prove that the limit surface carries indeed a vertical saddle connection (that belongs to the Delaunay polygonation). By definition of being Delaunay-convergent, all the Delaunay polygons of the sequence $g_{t_{n}}\Sigma$ converge toward a limit polygon of $\Sigma_{\infty}$. A polygon can only converge toward a polygon of equal or less number of sides so that the complexity can only decrease. Assuming that the complexity of $\Sigma_{\infty}$ and $\Sigma$ are the same, all the limit polygons keep the same number of sides. In particular, the side corresponding to the vertical connection does not vanish, and the limit surface carries a vertical saddle connection as well.
\end{proof}

On a dilation surface, vertical saddle connections have a top and a bottom endpoint. For any vertical saddle connection $\gamma$, we denote by $\mathcal{R}(\gamma)$ the vertical ray satisfying the following properties:
\begin{itemize}
    \item the starting point $M$ of $\mathcal{R}(\gamma)$ is the top endpoint of saddle connection $\gamma$;
    \item at $M$, $\mathcal{R}(\gamma)$ and $\gamma$ form an angular sector of amplitude $\pi$ contained in the right half-plane (by convention). In other words, $\mathcal{R}(\gamma)$ is obtained from $\gamma$ by turning counter-clockwise around $M$ by an angle of $\pi$.
\end{itemize}

We will prove that $\Sigma_{\infty}$ contains a cylinder with vertical boundary saddle connection by exhibiting a cyclic sequence of vertical saddle connections.

\begin{lem}\label{lem.demdivergedichotomy2}
For any vertical saddle connection $\gamma$ in the limit surface $\Sigma_{\infty}$, $\mathcal{R}(\gamma)$ is a vertical saddle connection too.
\end{lem}

\begin{proof}
We argue by contradiction, assuming that $\Sigma_{\infty}$ contains a vertical saddle connection $\gamma$ such that $\mathcal{R}(\gamma)$ is not a saddle connection. We denote by $s_{\mathrm{top}}$ the top singularity of $\gamma$ and choose a continuous parametrization $\Gamma(t)$ of the ray $\mathcal{R}(\gamma)$ with $\Gamma(0) = s_{\mathrm{top}}$, see Figure \ref{makeitdegsmin}. For $t$ small enough we define $\mathcal{D}_t$ the unique immersed disk whose center is $\Gamma(t)$ and whose radius is given by the segment $[s_{\mathrm{top}}, \Gamma(t)]$. Note that for any $t > t'$ we have $\mathcal{D}_{t'}  \subset \mathcal{D}_t$. We then define the "closest" singularity $s_{\mathrm{top}}^{\min}$ to $\gamma_{\infty}$ as the first singularity encountered when considering the increasing sequence of disks $(\mathcal{D}_t)_{t \geq 0}$. Note that such a sequence must actually encounter a singularity because of Theorem~\ref{thm:Veech}.\\

By hypothesis, $s_{\mathrm{top}}^{\min}$ does not belong to $\mathcal{R}(\gamma)$. We will reach our contradiction by showing that there is a closer singularity to $s_{\mathrm{top}}$ than $s_{\mathrm{top}}^{\min}$. In order to do so, we rely on the assumption that $\Sigma_{\infty}$ is the Delaunay limit of $(g_{t_{n}}\Sigma)_{n \in \mathbb{N}}$. All polygons of the Delaunay triangulation converge toward their limit polygon, all the quantities indexed by $n$ must converge toward their $\infty$-indexed corresponding quantity. For $N>0$ large enough, any surface $g_{t_{n}}\Sigma$ contains a well-defined saddle connection $\gamma_{n}$ corresponding to $\gamma_{\infty}$ in $\Sigma_{\infty}$. Analogously we denote by $s_{\mathrm{top}}(n)$ and $s^{\min}_{\mathrm{top}}(n)$ the top singularity of such a sequence and the singularity in $g_{t_{n}}\Sigma$ corresponding to $s^{\min}_{\mathrm{top}}$, see Figure \ref{makeitdegsmin}.  \\

\begin{figure}[h!]
\begin{center}
	\def\svgwidth{1 \columnwidth}
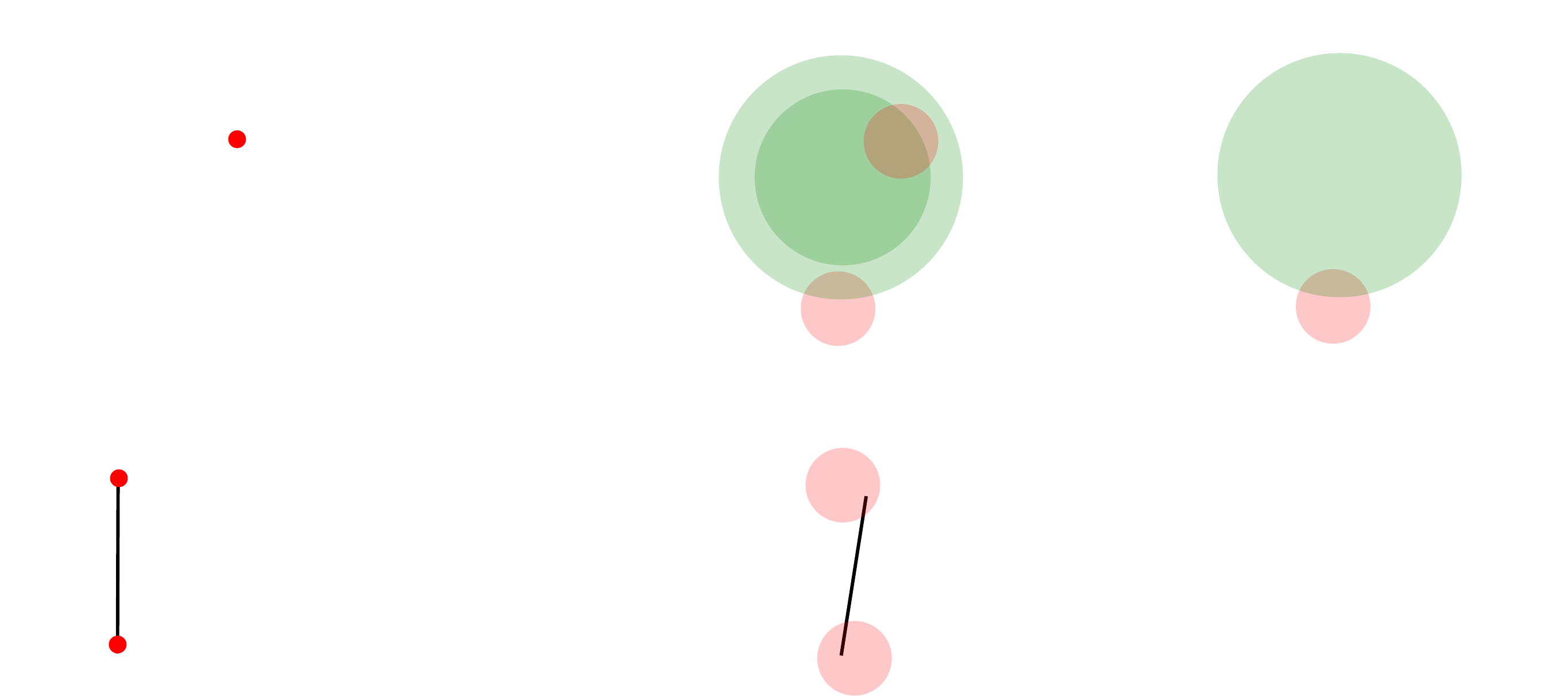
	\end{center}
\caption{On the left. In green, the maximal disk defining the singularity $s_{\mathrm{top}}^{\min}$. The red dots correspond to the singularities encountered on a $\mathcal{R}$-orbit in $\Sigma_{\infty}$. The blue dots correspond to the associated singularities in $\Sigma_n$. On the right. the left part corresponds to the surface $\Sigma_n$ as the right one to $\Sigma_m$. The shadow of the closest singularity on $\Sigma_m$ is mapped within a $2 \epsilon$ shorter disk of $\Sigma_n$.} 
	\label{makeitdegsmin}
\end{figure}

Let then $m > n$ such that $t_{m} - t_{n} > 0 $. By construction we have $g_{t_{n} - t_{m}}\Sigma_m = \Sigma_n$. Note that $t_{n} - t_{m} < 0 $ so that the Teichm\"{u}ller flow is now expanding (by a definite amount independent on $n$, $m$ and $\epsilon$) in the vertical direction and contracting the horizontal one. Note also that the image of $s^{\min}_{\mathrm{top}}(m)$ under $g_{t_{n} - t_{m}}$ must be a singularity of $\Sigma_n$. Since $\epsilon$ is arbitrary, one can take it as small as for the singularity $g_{t_{n} - t_{m}}(s^{\min}_{\mathrm{top}}(m))$ to be inside the disk centred at the same point that the maximal disk defining $s^{\min}_{\mathrm{top}}$ but of radius $2\epsilon$ shorter. This contradicts our initial assumption: as $g_{t_{n} - t_{m}}(s^{\min}_{\mathrm{top}}(m))$ must be $\epsilon$ close to a singularity of $\Sigma_{\infty}$, this new singularity would be inside the maximal disk defining $s^{\min}_{\mathrm{top}}$, see Figure~\ref{makeitdegsmin}.
\end{proof}

We now are able prove the main result of this Section.

\begin{proof}[Proof of Proposition~\ref{prop.degenerating}]
Assuming that $\Sigma_{\infty}$ and $\Sigma$ have the same complexity, Lemma~\ref{lem.demdivergedichotomy} proves that $\Sigma_{\infty}$ contains a vertical saddle connection $\gamma_{\infty}$. Using repeatedly Lemma~\ref{lem.demdivergedichotomy2}, we get that $\gamma_{\infty}$ belongs to a periodic sequence of vertical saddle connections. In other words, two consecutive saddle connections in this sequence differ by an angle of $\pi$. The curve formed by these saddle connection becomes simple if moved slightly toward the right. Such a simple curve must be the boundary of a cylinder. Consequently, the vertical saddle connection $\gamma_{\infty}$ belongs to the boundary of some cylinder of $\Sigma_{\infty}$. Proposition~\ref{prop:cylinderOPEN} guarantees that for any $\epsilon > 0$ there is $N>0$ such that for any $n \geq N$, $g_{t_{n}}\Sigma$ contains a closed geodesic whose direction belongs to $]\frac{\pi}{2}-\epsilon,\frac{\pi}{2}+\epsilon[$. For positive times, the Teichm\"{u}ller flow $g_{t}$ expands the horizontal direction and shrinks the vertical direction. Thus, \textit{a fortiori}, the same holds for $\Sigma$.
\end{proof}

\section{Cyclic maximal domains of type 1 and 2 (Proof of Proposition \ref{prop.maximaldomains})}
\label{sec:control}

We will show that the existence of a cyclic maximal domain of type 1 or 2 (see Sections~\ref{sec:domain1} and~\ref{sec:chain2}) in the Delaunay limit $\Sigma_{\infty}$ of a Delaunay-convergent subsequence $(g_{t_{n}}\Sigma)_{n\in\mathbb{N}}$ of positive Teichm\"{u}ller orbit of a dilation surface $\Sigma$ implies the existence of closed geodesic in $\Sigma$ whose direction is arbitrarily close to the vertical direction.\\

We will prove Proposition~\ref{prop.maximaldomains} by contradiction. We first give estimates on the moduli and directions of cylinders in surfaces of the positive Teichm\"{u}ller orbit of a dilation surface without closed geodesics in a neighborhood of the vertical direction.

\begin{lem}\label{lem:ControlTeich}
For $\delta>0$, we consider a closed dilation surface $\Sigma$ that does not contain any closed geodesic whose direction belongs to interval  $]\frac{\pi}{2} - \epsilon, \frac{\pi}{2} + \epsilon[$.
\par
Then, there is a positive constant $C_{\delta} > 0$ such that for any $t \geq 0$ the modulus of any cylinder of the surface $g_{t}\Sigma$ is bounded above by $C_{\delta}$.
\par
Besides, for any $\epsilon>0$, there is a time $T_{\epsilon}$ such that for any $t \geq T_{\epsilon}$, directions of closed geodesics of $g_{t}\Sigma$ are contained in $[\pi-\epsilon,\pi]$.
\end{lem}

\begin{proof}
The second claim follows immediately from the action of the Teichm\"{u}ller on interval $]\frac{\pi}{2} - \epsilon, \frac{\pi}{2} + \epsilon[$ in $\mathbb{RP}^{1}$.\\

Recall that Corollary~\ref{cor:moduli} asserts that every cylinder of $\Sigma$ is of modulus at most $M$ for some $M>0$. We will prove that the moduli of cylinders of surfaces $(g_{t}\Sigma)_{t\in\mathbb{R}^{+}}$ satisfy a global upper bound $C_{\delta}=\frac{M}{\sin^{2}(\delta)}$. Note that cylinders of $g_{t}\Sigma$ are correspond to cylinders of $\Sigma$.\\

We first consider a flat cylinder $\mathcal{C}$ of $\Sigma$. Normalising its area to $1$, its modulus is equal to $h^{-2}$ where $h$ is the (normalised) length of its closed geodesics. These geodesics have a direction $\theta$ which is $\delta$ far away from $\frac{\pi}{2}$ by assumption. Now we consider the images of $\mathcal{C}$ under the action of Teichm\"{u}ller flow. The normalised area remains identical while the lengths $h_{t}$ of closed geodesics of $g_{t}\mathcal{C}$ satisfy 

$$\frac{h_{t}}{h}=\sqrt{e^{2t}\cos^{2}(\theta)+e^{-2t}\sin^{2}(\theta)} \ . $$

It follows that 
    $$\frac{1}{h_t} \leq \frac{1}{\sin^2(\delta) \ h} \ ,$$
which is desired result. \\

It then remains to deal with the case where the cylinder $\mathcal{C}$ is a dilation cylinder of $\Sigma$ whose directions of closed geodesics is the interval $]\theta_{1},\theta_{2}[ \subset ]-\frac{\pi}{2} + \delta, \frac{\pi}{2} - \delta[$. We denote by $\lambda>1$ the dilation multiplier of $\mathcal{C}$. Recall that the modulus $M$ of such a cylinder is given by the relation 
$$ M = \frac{|\theta_{2}-\theta_{1}|}{\ln(\lambda)}  \ . $$
The action of Teichm\"{u}ller flow preserves the dilation multiplier. On the other hand, $g_{t}$ transforms any slope $\theta$ into 
$$\arctan(e^{-2t} \tan(\theta)) \ .$$

Using that for any $x,y <0$, we have $|\arctan(y)-\arctan(x)|\leq |y-x|$, we get that the size $|\theta_1(t) - \theta_2(t)|$ of the interval $g_{t}(](\theta_{1},\theta_{2}[)$ satisfies 
\begin{align*}
    |\theta_1(t) - \theta_2(t)| & \leq e^{-2t}| \tan(\theta_{1})-\tan(\theta_{2})| \\
    & \leq |\theta_{1}- \theta_{2}| \supr{\theta \in]-\frac{\pi}{2} + \delta, \frac{\pi}{2} - \delta[} | \tan'(\theta) |. 
\end{align*} 
Actually, $\supr{\theta \in]-\frac{\pi}{2} + \delta, \frac{\pi}{2} - \delta[} | \tan'(\theta) | = \frac{1}{\sin^{2}(\delta)}$ so we have
$$|\theta_{1}(t) - \theta_{2}(t)| \leq \frac{|\theta_{1}- \theta_{2}|}{\sin^{2}(\delta)}.$$
Thus, the modulus of the image cylinder is bounded above by $\frac{M}{\sin^{2}(\delta)}$.
\end{proof}

We split the proof of Proposition~\ref{prop.maximaldomains} in two statements, corresponding to the maximal domains of type 1 and type 2 respectively.

\subsection{Cyclic maximal domains of type 1}

\begin{figure}[h!]
\includegraphics[scale=0.6]{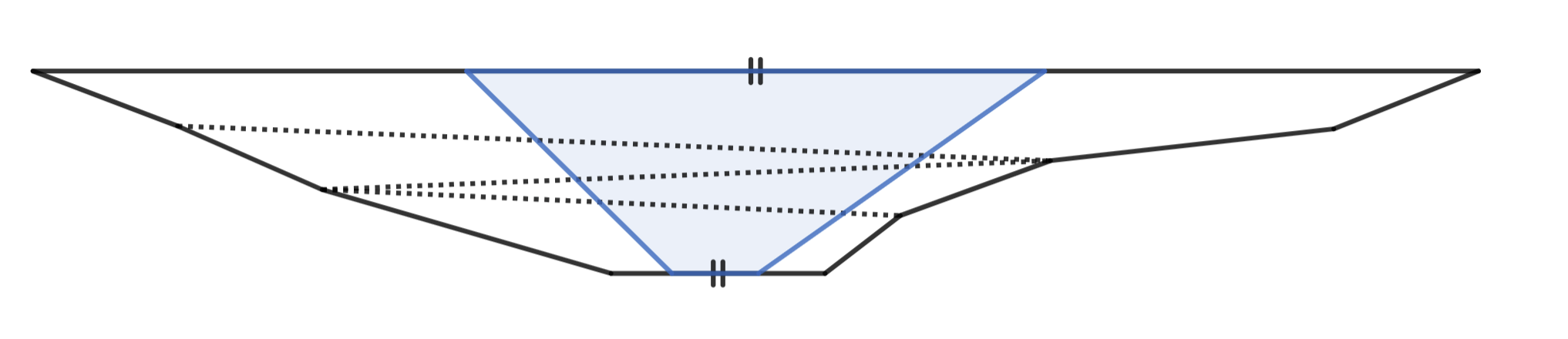}
\caption{A maximal domain of type 1 containing a dilation cylinder.}
\label{figure13}
\end{figure}

Polygons of type 1 assemble into maximal domains of type 1 (see Section~\ref{sec:domain1}). Following Proposition~\ref{prop:longside}, the (unique) long side of a degenerating polygon of type 1 must be glued to the short side of any other polygon. A cyclic maximal domain is formed by polygons of type 1 glued long side on short side, as in Figure~\ref{figure13}.

We prove that cyclic maximal domains of type 1 contain dilation cylinders whose angular amplitude is bounded by below.

\begin{prop}
\label{prop.maximaldomain1}
Let $\Sigma$ be a dilation surface and $t_n \to +\infty$ such that $(g_{t_n}\Sigma)_{n \in \mathbb{N}}$ Delaunay-converges and such that at least one of its Delaunay pieces is a cyclic maximal domain of type 1. Then for any $\epsilon> 0$, $\Sigma$ carries a cylinder whose direction belongs to $]\frac{\pi}{2} - \epsilon, \frac{\pi}{2} + \epsilon[$.
\end{prop}

\begin{proof}
Assuming by contradiction that for some $\epsilon>0$, interval $]\frac{\pi}{2} - \epsilon, \frac{\pi}{2} + \epsilon[$ does not contain any direction of closed geodesic of $\Sigma$, we use Lemma~\ref{lem:ControlTeich} to prove that the maximal angular amplitude of dilation cylinders of $\Sigma_{n}=g_{t_n}\Sigma$ becomes arbitrarily small as $n$ tends to infinity. We will obtain a contradiction by proving that
the cyclic maximal domain of type 2 $(X_{n})_{n \in \mathbb{N}}$ in $(\Sigma_{n})_{n\geq \mathbb{N}}$ contains a dilation cylinder whose angular amplitude is bounded by below provided $n$ is large enough.\\

The incidence graph of $(X_{n})_{n \in \mathbb{N}}$ is connected and contains a unique (oriented) cycle $C$ (see Section~\ref{sec:domain1} for details). For any $n \in \mathbb{N}$, $X_{n}$ is a topological cylinder. We consider an edge $(L_{n})_{n \in \mathbb{N}}$ between two polygons of the cycle $C$.\\

Cutting along edge $(L_{n})_{n \in \mathbb{N}}$ in $(X_{n})_{n \in \mathbb{N}}$, we obtain a sequence of simply connected flat surfaces $(P_{n})_{n \in \mathbb{N}}$ with a unique boundary component. It is formed by the gluing of polygons of type 1 according to an incidence graph which is a tree.\\

By definition of a polygon of type 1 (see Definition~\ref{defn:DCpolygons}), all the sides of $(P_{n})_{n \in \mathbb{N}}$ have the same limit direction in $\mathbb{RP}^{1}$. We normalize $(P_{n})_{n \in \mathbb{N}}$ in such a way that all the sides tend to be horizontal and every surface $P_{n}$ has unit area. In particular, provided $n$ is large enough, $P_{n}$ is a planar polygon in the classical sense.\\

Since every polygon of type 1 has a unique long side (see Definition~\ref{defn:longshort}), provided that $n$ is large enough, $P_{n}$ has a unique upper side $S_{n}$ (or a unique lower side, depending on the normalization) and several lower sides $T_{n}^{1},\dots,T_{n}^{p-1}$ (where $p$ is the number of sides of $P_{n}$ for any $n$).\\

Edge $(L_{n})_{n \in \mathbb{N}}$ corresponds to the identification of the unique upper side $(S_{n})_{n \in \mathbb{N}}$ with some lower side $(T_{n}^{i_{0}})_{n \in \mathbb{N}}$. Since $S_{n}$ and $T_{n}^{i_{0}}$ have the same slope, they cannot be adjacent in the boundary of $P_{n}$. Thus, the corner angles and the ends of $T_{n}^{i_{0}}$ tend to $\pi$ as $n$ tends to infinity. Provided that $n$ is large enough, rays starting from the ends of side $T_{n}^{i_{0}}$ in directions $\frac{\pi}{4}$ and $\frac{3\pi}{4}$ intersect $S_{n}$ and there exists a trapezoid $M_{n}$ in $P_{n}$ formed by $T_{n}^{i_{0}}$ (the lower side of $M_{n}$), a side of slope $\frac{\pi}{4}$, a portion of $S_{n}$ (the upper side of $M_{n}$) and a side of slope $\frac{3\pi}{4}$.\\

Since any point of the upper side of $M_{n}$ is identified with a point of $T_{n}^{i_{0}}$, the lift of trapezoid $M_{n}$ in $X_{n}$ contains a family of closed geodesics whose slopes sweep an interval of length at least $\frac{\pi}{2}$ in $\mathbb{RP}^{1}$ (see Figure~\ref{figure13}). Thus, for any large enough $n$, surface $\Sigma_{n}$ contains a dilation cylinder of angle at least $\frac{\pi}{2}$. This is the desired contradiction.
\end{proof}

\subsection{Cyclic maximal domains of type 2}

\begin{figure}[h!]
\includegraphics[scale=0.6]{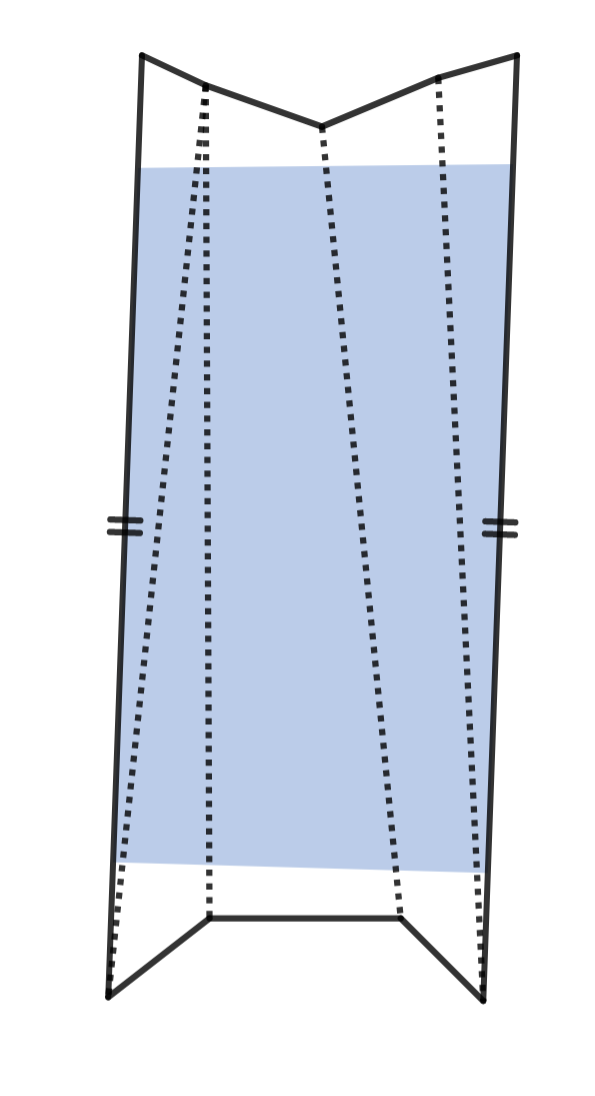}
\caption{A maximal domain of type 2 containing a cylinder of large modulus (in blue).}
\label{figure14}
\end{figure}

In Section~\ref{sec:chain2}, we defined maximal domains of type 2 as collections of polygons of type 2 glued along their long boundary sides. Such a maximal domain is cyclic if the polygons are glued according to a cyclic graph, as in Figure~\ref{figure14}.

We prove that cyclic maximal domains of type 2 contain cylinders of arbitrarily large modulus.

\begin{prop}
\label{prop.maximaldomain2}
Let $\Sigma$ be a dilation surface and $t_n \to +\infty$ such that $(g_{t_n}\Sigma)_{n \in \mathbb{N}}$ Delaunay-converges and such that at least one of its Delaunay pieces is a cyclic maximal domain of type 2. Then for any $\epsilon> 0$, $\Sigma$ carries a cylinder whose direction belongs to $]\frac{\pi}{2} - \epsilon, \frac{\pi}{2} + \epsilon[$.
\end{prop}

\begin{proof}
We proceed similarly as when proving Proposition~\ref{prop.maximaldomain1}. Assuming by contradiction that for some $\epsilon>0$, interval $]\frac{\pi}{2} - \epsilon, \frac{\pi}{2} + \epsilon[$ does not contain any direction of closed geodesic of $\Sigma$, we use Lemma~\ref{lem:ControlTeich} to obtain an upper bound $M>0$ on the modulus of any cylinder in any dilation surface $\Sigma_{n}=g_{t_n}\Sigma$. We will prove that the cyclic maximal domain of type 2 $(X_{n})_{n \in \mathbb{N}}$ in $(\Sigma_{n})_{n \in \mathbb{N}}$ contains cylinders of arbitrarily large modulus as $n$ tends to infinity.\\

In particular, for any $n \in \mathbb{N}$, $X_{n}$ is a topological cylinder. We cut along some edge $(L_{n})_{n \in \mathbb{N}}$ and obtain a sequence of polygons $(P_{n})_{n \in \mathbb{N}}$. We normalize each polygon $P_{n}$ in such a way that the two sides corresponding to edge $L_{n}$ are vertical and $P_{n}$ has unit area. These two vertical sides will be referred to as $S_{n}$ (for the left side) and $T_{n}$ (for the right side).\\

In a polygon of type 2 (see Definition~\ref{defn:DCpolygons}), the ratio between the length $|S_{n}|$ of $S_{n}$ and the length $|T_{n}|$ of $T_{n}$ converges to $1$ and the length of any other side of $P_{n}$ becomes negligible in comparison with $|S_{n}|$ and $|T_{n}|$ (see Figure~\ref{figure14}). Since $P_{n}$ has unit area for any $n \in \mathbb{N}$, $|S_{n}|$ and $|T_{n}|$ tends to infinity while the distance between $S_{n}$ and $T_{n}$ tends to zero as $n \rightarrow + \infty$.\\

It follows that for any $\epsilon>0$, there is $N>0$ such that for any $n \geq N$, polygon $P_{n}$ contains a rectangle $R_{n}$ satisfying the following conditions:
\begin{itemize}
    \item sides of $R_{n}$ are either vertical or horizontal;
    \item the vertical left and right sides are portions of $S_{n}$ and $T_{n}$;
    \item the length of the vertical sides of $R_{n}$ is at least $(1-\epsilon)|S_{n}|$;
    \item the length of the vertical sides of $R_{n}$ tend to infinity as $n \rightarrow + \infty$;
    \item the length of the horizontal sides of $R_{n}$ tend to zero as $n \rightarrow + \infty$.
\end{itemize}
Since sides $S_{n}$ and $T_{n}$ are identified, the lift of rectangle $R_{n}$ in $X_{n}$ contains a family of closed geodesics covering most of rectangle $R_{n}$ (the complement of a part of arbitrarily small relative area). Therefore, provided that $n$ is large enough, $X_{n}$ contains a cylinder of arbitrarily large modulus (see Figure~\ref{figure14}). We obtain the desired contradiction.
\end{proof}

\section{Long sides and short sides (Proof of Proposition~\ref{prop.petiteportegrandeporte})}\label{sec.petiteportegrandeporte}

We will actually prove the following slightly stronger version of Proposition
\ref{prop.petiteportegrandeporte}, which does not involve the Teichmüller flow. \\

\textbf{Proposition 4.4}. \textit{Let $(\Sigma_n)_{n \in \mathbb{N}}$ be a Delaunay-convergent sequence of dilation surfaces such that all its Delaunay pieces have at least one long boundary side. Then, for any open set $U \subset \mathbb{RP}^{1}$, there is $N>0$ such that for any $n \geq N$, $\Sigma_n$ contains closed geodesics whose direction belong to $U$.}  \\

The proof is based on Proposition~\ref{prop:boundaryexceptional}. This proposition asserts that, in the case of a dilation surface with boundary, either a given open set of directions contains a cylinder or a set of trajectories having these directions, a pencil to be precise, must leave across a boundary component of the dilation surface. This proposition then shows that we can concentrate on the case where each trajectory of a Delaunay piece (see Definition~\ref{defn:PIECES}) leaves it by hitting the boundary. It can do it by crossing either a long edge or a short one. We will actually rule out the short edge case in Sections~\ref{sub:traj12} and~\ref{sub:trajcore}. Indeed, these boundaries are by definition very small compared to the long edges and it will be unlucky to leave the piece through such a short side. The trajectories of the pencil will then have to leave the Delaunay piece through a long side and then enter a new Delaunay piece through a short side (as by construction Delaunay pieces are glued to one another short side to long side, see Lemma~\ref{lem:Delpiece}). If the pencil does not enter in a cylinder one can repeat the argument to get a sequence of Delaunay pieces such that the pencil enters them by short side and leave them by long sides. As there are only finitely many boundary components, such a pencil will cross twice a given edge. The first return map on such an edge is a very dilating mapping as going from short sides to long sides induces a huge contraction. This concludes as contracting mappings have periodic orbits.

\subsection{Trajectories inside maximal domains of type 1 or 2}\label{sub:traj12}

For a trajectory whose slope is far enough from the limit directions of the (finitely many) Delaunay edges, we have some control on its behaviour in Delaunay pieces formed by degenerating polygons.

\begin{defn}\label{defn:TRANSVERSE}
For any $\epsilon>0$, $\Theta_{\epsilon} \subset \mathbb{RP}^{1}$ is the open subset of slopes whose distance to any limit direction of a Delaunay edge of $(\Sigma_{n})_{n \in \mathbb{N}}$ is strictly bigger than $\epsilon$.
\end{defn}

The following proposition asserts that a for a given direction in $\Theta_{\epsilon}$ a trajectory entering a maximal domain of type 1 by a small edge must exit it through a long one, provided that $n$ is large enough.

\begin{prop}\label{prop:transverse1}
Let $(E_{n})_{n \in \mathbb{N}}$ be a short boundary edge of a maximal domain of type~1 $(X_{n})_{n \in \mathbb{N}}$ that has at least one boundary long edge.
\par
For any $\epsilon \in ]0,\frac{\pi}{2}[$, there is a long boundary edge $(M_{n})_{n \in \mathbb{N}}$ of $(X_{n})_{n \in \mathbb{N}}$ and $N>0$ such that for any $n \geq N$, any trajectory of $X_{n}$ whose direction belongs to $\Theta_{\epsilon}$ starting from $E_{n}$ eventually leaves $X_{n}$ through the interior of $M_{n}$.
\end{prop}

\begin{proof}
Since a maximal domain of type 1 is formed by polygons of type 1, Delaunay edges have the same limit slope. Without loss of generality, we will assume that this unique limit slope is horizontal.
\par
We start by discussing the case of a polygon of type 1. Note that for any $\delta > 0$, there is $N_{\delta} \in \NN$ such that for any $n \geq N_{\delta}$, every (convex) Delaunay polygon $P_{n}$ of $X_{n}$ satisfies the following properties:
\begin{itemize}
    \item the slope of every Delaunay edge belongs to $]-\delta,\delta[ \subset \mathbb{RP}^{1}$;
    \item the inner angle between two short sides of $P_{n}$ is at least $\pi - \delta$;
    \item the inner angle between a short side and a long side of $P_{n}$ is at most $\delta$.
\end{itemize}
If a trajectory of $P_{n}$ starts from a short side and leaves $P_{n}$ through another short side, then it cuts out $P_{n}$ into two polygons. Computing the sum of the inner angles in each of them, we deduce that the slope of $t$ belongs to $[-p\delta,p\delta]$ where $p$ is the number of sides of $P_{n}$. Thus, by choosing $\delta \geq \frac{\epsilon}{q}$ where $q$ is the number of Delaunay edges of $P_{n}$, one makes sure that for $n \geq N_{\delta}$, a trajectory of polygon $P_{n}$ whose slope belongs to $\Theta_{\epsilon}$ starting from a short side of $P_{n}$ leaves it through the interior of its unique long side. \\

The proof of the non cyclic domain of type 1 follows the exact same line. The key remark being that type 1 polygons piled up long side to short side form a polygon that satisfies the three points above (see Figure \ref{figure13}). Therefore, if we set $\delta = \frac{\epsilon}{m}$ where $m$ is the total number of edges that are short sides of at least one polygon of $X_n$, following the argumentation above, we see that there is $N_{\delta}$ large enough such that for $n \geq N_{\delta}$ the trajectory visits finitely many long boundaries of $X_n$ and exits $X_n$ as otherwise the domain would be cyclic.
\end{proof}

We now address the case of maximal domains of type 2.

\begin{prop}\label{prop:transverse2}
Let $(E_{n})_{n \in \mathbb{N}}$ be a short boundary edge of a maximal domain of type~2 $(X_{n})_{n \in \mathbb{N}}$ that has at least one boundary long edge. We also consider a non-empty open interval $\mathcal{I} \subset \Theta_{\epsilon}$ for some $\epsilon>0$.
\par
There is a long boundary edge $(M_{n})_{n \in \mathbb{N}}$ of $(X_{n})_{n \in \mathbb{N}}$ and $N>0$ such that for any $n \geq N$, any trajectory of $X_{n}$ whose direction belongs to $\mathcal{I}$ starting from $E_{n}$ eventually leaves $X_{n}$ through the interior of $M_{n}$.
\end{prop}

\begin{proof}
A maximal domain of type 2 with at least one boundary edge actually has two long boundary edges since its graph of incidence is linear (see Section~\ref{sec:chain2}). Without loss of generality, we assume that the limit slope of the long Delaunay edges of $(X_{n})_{n \in \mathbb{N}}$ is vertical. Therefore, we will refer to the boundary edges of $X_{n}$ (which is a polygon) as the long left side, the long right side, the short upper sides and the short lower sides. \\

There is $N>0$ such that for any $n \geq N$, the slope of any straight segment joining an upper vertex and a lower vertex of $X_{n}$ is contained in $]\frac{\pi}{2}-\epsilon,\frac{\pi]}{2}+\epsilon[$. It follows that the slope of a trajectory of $X_{n}$ joining a short upper side and a short lower side cannot belong to $\Theta_{\epsilon}$ for $n \geq N$. It remains to consider the case of a trajectory $t$ joining two upper (or lower) short sides of $X_{n}$. \\

It follows from Proposition~\ref{prop:concave} that for any $\delta>0$, there is $N_{\delta}$ such that for any $n \geq N_{\delta}$, the inner angle between two consecutive upper (or lower) short sides of $X_{n}$ is at least $\pi - \delta$. 
Any trajectory $t$ cuts out $X_{n}$ into two polygons. Computing the sum of inner angles in each of them, we deduce that trajectory $t$ form an angle of magnitude smaller than $p\delta$ with one of the short sides of $X_{n}$ (here, $p$ is the number of sides of $X_{n}$). Since $\delta$ can be made arbitrarily small, there exists a bound $N'>0$ such that for any $n \geq N'$, a trajectory joining two upper (or lower) short sides of $X_{n}$ cannot belong to $\Theta_{\epsilon}$.
\par
Consequently, for any $n$ satisfying $n \geq \max(N,N')$, any trajectory starting from a short boundary edge $E_{n}$ of $X_{n}$ eventually leaves $X_{n}$ through the interior of one of its two extremal edges (see Figure~\ref{figure15}). If we restrict ourselves to trajectories whose slope belongs to a connected open subset $U$ of $\Theta_{\epsilon}$, a continuity argument proves that two trajectories starting from $E_{n}$ leave $X_{n}$ through the same extremal edge.

\begin{figure}[h!]
\includegraphics[scale=0.6]{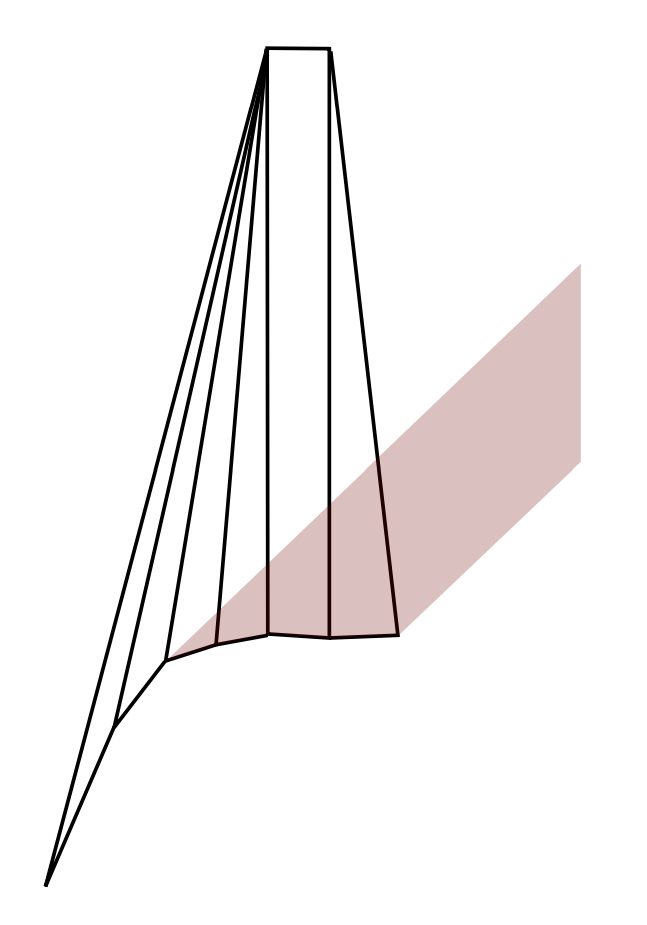}
\caption{A maximal domain of type 2 with trajectories starting from short sides in a direction far for that of the long sides.}
\label{figure15}
\end{figure}

\end{proof}

\subsection{Trajectories inside connected components of the core}\label{sub:trajcore}

The case of Delaunay pieces that are connected components of the core $(\mathcal{C}\Sigma_{n})_{n \in \mathbb{N}}$ is a bit more complicated. In order to find an open set of directions where trajectories starting from the same short boundary edge leave a component of $\mathcal{C}\Sigma_{n}$ through the same long boundary edge, we first prove the analogous result for connected components of the limit surface $\Sigma_{\infty}$, which is an easy consequence of Proposition~\ref{prop:boundaryexceptional}.

\begin{lem}\label{lem:leavingINFINITE}
For any non-empty open subset $U \subset \mathbb{RP}^{1}$ and any connected component $X_{\infty}$ of $\Sigma_{\infty}$ with a non-empty boundary, one of the following statements holds:
\begin{itemize}
    \item there exists a closed geodesic in $X_{\infty}$ whose slope is contained in $U$;
    \item there is a non-empty open subset $V \subset U$ such that every trajectory starting from a singularity $x$ of $X_{\infty}$ in a direction of $V$ eventually leaves $X_{\infty}$ through the interior of a boundary saddle connection.
\end{itemize}
\end{lem}

\begin{proof}
Let $S_{x,U}$ be the set of (oriented) trajectories starting from the singularity $x$ with a slope in $U$. The topology of $S_{x,U}$ is induced by the canonical projection $\pi_{x}$ to $\mathbb{RP}^{1}$. Assuming that no closed geodesic of $X_{\infty}$ belongs to a direction of $U$, it has been proved in Proposition~\ref{prop:boundaryexceptional} that trajectories of $S_{x,U}$ leaving $X_{\infty}$ through the interior of a boundary saddle connection form an open dense subset of $S_{x,U}$. Since there are finitely many such singularities in $\Sigma_{\infty}$ and projections $\pi_{x}$ have finitely many preimages, there is an open dense subset $V$ of $U$ such that every trajectory starting from such a singularity $x$ in a direction of $V$ leaves its component through the interior of a boundary saddle connection.
\end{proof}

Since short boundary edges of connected components of the core degenerate to singular points in the Delaunay limit, we obtain a result about trajectories in the connected components of the core.

\begin{prop}\label{prop:leaving}
Let $(X_{n})_{n \in \mathbb{N}}$ be a connected component of the core $(\mathcal{C}\Sigma_{n})_{n\in\mathbb{N}}$ with at least one long boundary edge. Let $(E_{n})_{n \in \mathbb{N}}$ be a short boundary edge of $(X_{n})_{n \in \mathbb{N}}$. For any non-empty open subset $U \subset \mathbb{RP}^{1}$, one of the following statements holds:
\begin{itemize}
    \item there is a bound $N>0$ such that for any $n \geq N$, there exists a closed geodesic in $X_{n}$ whose slope is contained in $U$;
    \item there is a non-empty open subset $V \subset U$, a bound $N>0$ and a long boundary edge $(M_{n})_{n \in \mathbb{N}}$ of $(X_{n})_{n \in \mathbb{N}}$ such that for any $n \geq N$, any trajectory of $X_{n}$ in a direction $\theta \in V$ starting from $E_{n}$ eventually leaves $X_{n}$ through the interior of $M_{n}$.
\end{itemize}
\end{prop}

\begin{proof}
We first decompose the proof into two subcases depending whether $\Sigma_{\infty}$ contains a closed geodesic whose direction belongs to $U$ or not. In the first case, we deduce from Proposition~\ref{prop:cylinderOPEN} that there exists $N>0$ such that for any $n \geq N$, $\Sigma_{n}$ contains a closed geodesic whose direction belongs to $U$. \\

In the second case, we fix $\epsilon$ small enough such that $U \cap \Theta_{\epsilon}$ is non-empty. Lemma~\ref{lem:leavingINFINITE} then proves the existence of a non-empty open subset $V$ of $U \cap \Theta_{\epsilon}$ such that any trajectory of the Delaunay limit $X_{\infty}$ of $(X_n)_{n \in \mathbb{N}}$ whose direction belongs to $V$ and starting from a singularity leaves $X_{\infty}$ through the interior of a boundary saddle connection. \\

Let $(E_n)_{n \in \mathbb{N}}$ be a short boundary edge of $(X_n)_{n \in \mathbb{N}}$. By construction this short boundary edge converges toward a point $x$ of $X_{\infty}$ in the limit. For any interval $I$ in $V$, we consider the two-parameters family $P(E_n,I)$ of trajectories starting from the edge $E_{n}$ and whose direction belongs to $I$. This family of trajectories accumulates on a pencil $P(x, I)$ of $X_{\infty}$ as $n$ tends to infinity. \\

The edge $(E_{n})_{n \in \mathbb{N}}$ is a short edge of a chain of polygons of type 2 belonging to $(X_{n})_{n \in \mathbb{N}}$ (see Figure~\ref{figure10}). Using Proposition~\ref{prop:concave} as in the proof of Proposition~\ref{prop.maximaldomain2}, we deduce that provided that $n$ is large enough, trajectories of $P(E_n,I)$ leave each of these polygons of type 2 through one of its long side (the hypothesis that $I$ is disjoint from $\Theta_{\epsilon}$ is crucial here). Then, these trajectories finally enter a polygon of type 3 of $(X_n)_{n \in \mathbb{N}}$. A continuity argument proves that trajectories of pencil $P(x, I)$ leave $X_{\infty}$ through the interior of the same boundary saddle connection $M_{\infty}$ which is the limit of a long boundary edge $(M_{n})_{n \in \mathbb{N}}$ of $(X_{n})_{n \in \mathbb{N}}$. Up to replacing $I$ by a smaller open interval, we can assume that the intersection of trajectories of $P(x, I)$ with $M_{\infty}$ is disjoint from a neighborhood of the endpoints of $M_{\infty}$. We deduce that provided $n$ is large enough, trajectories of $P(E_n, I)$ leave $X_{n}$ through the interior of $M_{n}$.
\end{proof}

We combine the previous results to exhibit a set of directions and a lower bound that hold for every Delaunay piece of $(\Sigma_{n})_{n \in \mathbb{N}}$.

\begin{cor}\label{cor:leaving}
For any non-empty open subset $U \subset \mathbb{RP}^{1}$, one of the following statements holds:
\begin{itemize}
    \item there is a bound $N>0$ such that for any $n \geq N$, there exists a closed geodesic in $\Sigma_{n}$ whose slope is contained in $U$;
    \item there is a non-empty open subset $V \subset U$ and a bound $N>0$ such that for any $n \geq N$ and any short boundary edge $E_{n}$ in any Delaunay piece $X_{n}$ of $\Sigma_{n}$ having a long boundary edge, there is a long boundary edge $M_{n}$ such that every trajectory of $X_{n}$ starting from $E_{n}$ and whose slope belongs to $U$ eventually leaves $X_{n}$ through the interior of $M_{n}$.
\end{itemize}
\end{cor}

\begin{proof}
Provided $\epsilon$ is small enough, $\Theta_{\epsilon} \cap U$ is non-empty. For such a small $\epsilon$, we consider an open interval $I \subset \Theta_{\epsilon} \cap U$. Since there are finitely many Delaunay pieces and Delaunay edges in $(\Sigma_{n})_{n \in \mathbb{N}}$, there is a global bound $N_{0}$ such that the second statement holds for trajectories whose slope is in interval $I$ for any short boundary edge in any Delaunay piece $X_{n}$ that is a maximal domain of type 1 or 2 (see Propositions~\ref{prop:transverse1} and~\ref{prop:transverse2}). \\

Then, we apply Proposition~\ref{prop:leaving} to a boundary short edge $(E_{n})_{n \in \mathbb{N}}$ in a connected component $(X_{n})_{n \in \mathbb{N}}$ of the core. If $X_{n}$ contains a closed geodesic provided $n$ is large enough, then the first statement of our proposition holds. Otherwise, the second statement of Proposition~\ref{prop:leaving} provides a non-empty open subset $I'$ of $I$ and a new bound such that the property also holds for this edge. After finitely many steps, we obtain a nonempty open subset $V$ of $\mathbb{RP}^{1}$ and a bound $N>0$ such that the property holds of trajectories whose slope belongs to $V$ for every short boundary edge in every Delaunay piece having a long boundary edge (provided that $n \geq N$).
\end{proof}

\textbf{Proof of Proposition \ref{prop.petiteportegrandeporte}}. Corollary \ref{cor:leaving} shows that it is enough to prove that any direction $d$ such that for $n$ large enough any trajectory starting from any short edge of any Delaunay piece exits the Delaunay piece through one of its long sides carries a cylinder. As we only have finitely many Delaunay pieces, any such trajectory will have to cross twice some boundary edge $(E_n)_{n \in \mathbb{N}}$ of two Delaunay pieces (one for which it is a short side and one for which it is a long side). \\

\begin{figure}[h!]
\begin{center}
	\def\svgwidth{1 \columnwidth}
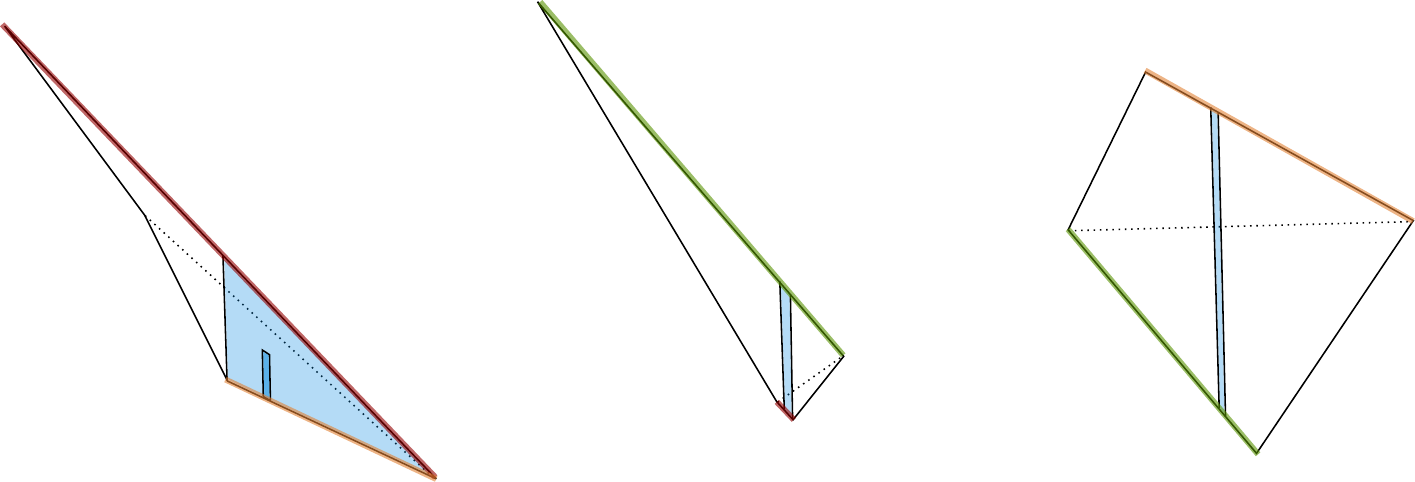
	\end{center}
\caption{From the left to the right: a maximal domain $D_1$ of type 1, a maximal domain $D_2$ of type 2 and a component $C$ of the core. The orange edge of the left corresponds to $E_n$. The ribbon of parallel trajectories leaves the $D_1$ to enter $D_2$ through a its short side which goes to $C$ from a long side of $D_2$. When the ribbon enters back $D_1$ it has been contracted by an amount that goes to infinity when $n \to +\infty$.}
	\label{fig.premierretour}
\end{figure}
Given a direction $d$ and a short edge $(E_n)_{n \in \mathbb{N}}$ we denote by $ P(E_n, d) $ the set of trajectories starting from a point of $(E_n)_{n \in \mathbb{N}}$ of direction $d$ pointing inside the Delaunay piece for which $(E_n)_{n \in \mathbb{N}}$ is the short side. \\

By definition of $E_n$, there is a trajectory $t$ of $P(E_n, d)$ which crosses back $E_n$ for $n$ large enough. We claim that all the trajectories of $P(E_n, d)$ cross $E_n$ alongside with $t$. Indeed, as any trajectory of $P(E_n, d)$ only exits a Delaunay piece by its long side and enters one by its short side, the contraction ratio of the first return map on a neighbourhood of $E_n$ converges to $0$ as $n \to + \infty$. In particular, for $n$ large enough, the image of $E_n$ must be fully contained in $E_n$, which implies that it has a periodic point. This periodic point of the first return map corresponds to a closed geodesic, concluding. \hfill $\qed$

\nopagebreak
\vskip.5cm

\end{document}